\documentclass{amsart}
\usepackage{latexsym, amsfonts, amscd, amsmath, verbatim, graphicx, amsthm}

\def\be{\begin{enumerate}}
\def\ee{\end{enumerate}}
\def\ii{\item}
\def\bi{\begin{itemize}}
\def\ei{\end{itemize}}

\def\RR{\mathbb{R}}

\def\beqa{\begin{eqnarray*}}
\def\eeqa{\end{eqnarray*}}
\def\beqan{\begin{eqnarray}}
\def\eeqan{\end{eqnarray}}

\def\ZZ{\mathbb{Z}}

\def\RR{\mathbb{R}}

\def\GG{\mathcal{G}}
\def\MM{\mathcal{M}}

\def\del{\partial}

\newtheorem{thm}{Theorem}[section]

\newtheorem{lemma}[thm]{Lemma}

\newtheorem{prop}[thm]{Proposition}
\newtheorem{cor}[thm]{Corollary}
\newtheorem{lemmadef}[thm]{Definition-Lemma}

\newtheorem{defn}[thm]{Definition}

\newtheorem{remark}[thm]{Remark}

\def\proof{\textit{Proof:}\ \ \ \ }

\def\qed{\hfill $\eop$ \vskip10pt}

\def\prob[#1]{\subsection*{-----\textsf{#1}-----}}
\def\sec[#1]{\noindent ---{\large \textit{#1}}--- \\}
\def\secv[#1]{\vspace{.3in} \noindent ---{\large \textit{#1}}--- \\}

\def\Hom{\textrm{Hom}}

\def\Re{\textrm{Re}}
\def\Im{\textrm{Im}}

\def\det{\textrm{det}}

\def\lim{\textrm{lim}}

\def\ker{\textrm{ker}}
\def\ind{\textrm{ind}}

\def\id{\textrm{id}}

\def\Symp{\textrm{Symp}}
\def\Diff{\textrm{Diff}}

\def\Nov{\textrm{Nov}}
\def\Flux{\textrm{Flux}}
\def\Ham{\textrm{Ham}}

\def\eop{\rule{5pt}{5pt}}

\title[Fixed points of surface symplectomorphisms]{A Sharp lower bound on fixed points of surface symplectomorphisms in each mapping class}

\author[A. Cotton-Clay]{Andrew Cotton-Clay}

\begin{document}

\begin{abstract}
Given a closed, oriented surface $\Sigma$, possibly with boundary, and a mapping class, we obtain sharp lower bounds on the number of fixed points of a surface symplectomorphism (i.e. area-preserving map) in the given mapping class, both with and without nondegeneracy assumptions on the fixed points. This generalizes the Poincar\'e-Birkhoff fixed point theorem to arbitrary surfaces and mapping classes. These bounds often exceed those for non-area-preserving maps. We obtain these bounds from Floer homology computations with certain twisted coefficients plus a method for obtaining fixed point bounds on entire symplectic mapping classes on monotone symplectic manifolds from such computations. For the case of possibly degenerate fixed points, we use quantum-cup-length-type arguments for certain cohomology operations we define on summands of the Floer homology.
\end{abstract}

\maketitle

\section{Introduction}

\subsection{Statement of Results}

Let $(\Sigma,\omega)$ be a compact surface of negative Euler characteristic\footnote{Our results extend to the nonnegative Euler characteristic case but these exceptional cases would be cumbersome to carry around. All of these cases are already understood.}, possibly with boundary, with $\omega$ a symplectic form (i.e. area form). For any mapping class we give sharp lower bounds on the number of fixed points of an area-preserving map $\phi$ in the mapping class, both in the case in which $\phi$ is assumed to have nondegenerate fixed points\footnote{That is, the fixed points of $\phi$ are cut out transversally in the sense that $\det(1-d\phi_x) \neq 0$ for fixed points $x$.} and in the general case in which degenerate fixed points are allowed. This generalizes the Poincar\'e-Birkhoff fixed point theorem, which states that area-preserving twist maps of the annulus have at least two fixed points, to arbitrary surfaces and mapping classes.

A traditional lower bound, which is sharp for non-area-preserving maps on surfaces, comes from Nielsen theory. The Nielsen class $\eta \in \pi_0(\Gamma(M_\phi))$ of a fixed point $x$ of a map $\phi$ is its homotopy class when thought of as a constant section of the mapping torus $M_\phi \rightarrow S^1$. The index of a Nielsen class $\eta$, denoted $\ind(\eta)$, is given by the sum of the (differential topology) indices of each fixed point in the class $\eta$. This quantity is invariant under deformation. The traditional lower bound on fixed points when there is a nondegeneracy condition is given by $\sum_\eta \left| \ind(\eta) \right|$, as nondegenerate fixed points have index $\pm 1$. In the general case with no nondegeneracy condition, the traditional lower bound is given by the number of Nielsen classes with nonzero index.

For closed surfaces, mapping classes are elements of $\pi_0(\Diff^+(\Sigma))$. For surfaces with boundary, we consider $\Diff_\del(\Sigma)$, diffeomorphisms with no fixed points on the boundary, and use the term mapping classes to refer to elements of $\pi_0(\Diff_\del(\Sigma))$, though this is not standard. By Moser's trick, in dimension two, these are homotopy equivalent to the versions with $\Diff$ replaced by $\Diff_{\mathrm{Vol}} = \Symp$.

We state our results in terms of Thurston's classification of surface diffeomorphisms \cite{t}, \cite{flp}, which states that every element of $\pi_0(\Diff^+(\Sigma))$ is precisely one of
\bi
\ii Periodic (finite order): For some representative $\phi$, we have $\phi^\ell = \id$ for some $\ell \in \ZZ_{>0}$
\ii Pseudo-Anosov: Some singular representative $\phi$ preserves two transverse singular measured foliations, expanding the measure on one and contracting the measure on the other. See \cite{c} for symplectic smoothings of these singular representatives.
\ii Reducible: Some representative $\phi$ fixes setwise a collection of curves $C$ none of which are nulhomotopic or boundary parallel (and the mapping class is not periodic)
\ei
In the reducible case, cutting along a maximal collection of curves $C$ gives a map on each component of $\Sigma\setminus C$ (given by the smallest power of $\phi$ which maps that component to itself) which is periodic or pseudo-Anosov. We call these components of $\Sigma\setminus C$ \emph{reducible components}.

When $\Sigma$ has boundary, our definition of mapping classes on surfaces with boundary distinguishes additionally an affine space over $\ZZ^b$ for each element of $\pi_0(\Diff^+(\Sigma))$ with $b$ setwise fixed boundary components. The action of the coordinate generators of $\ZZ^b$ is by full Dehn twists at the appropriate boundary component. We thus consider each element of $\pi_0(\Diff^+(\Sigma))$ to come with additional data of which direction each setwise fixed boundary component rotates and some nonnegative number of annular reducible components (one for each Dehn twist in the direction given).

\begin{thm}
\label{fixedbound}
The minimum number of fixed points of an area-preserving map $\phi$ with nondegenerate fixed points in a mapping class $h$ is given by:
$$ \left\{ \begin{array}{cc}
\sum_\eta \left| \ind(\eta) \right| & h\ \textrm{is periodic or pseudo-Anosov} \\
\sum_\eta \left| \ind(\eta) \right| + 2 A & h\ \textrm{is reducible} \end{array} \right.$$
where $A$ is the number of genus zero components of the reducible mapping class on which the map is the identity, which do not abut any pseudo-Anosov components, and all of whose boundary components rotate in the same direction and are nulhomologous or homologically boundary parallel.
\end{thm}

The upper bound here is given by construction, which comes from perturbing maps which are nice with respect to the Nielsen-Thurston geometry, which we call standard form maps, with particular symplectic vector fields. The lower bound comes from Floer homology computations for these standard form maps with certain twisted coefficients, based on computations we performed in \cite{c}, plus a result discussed below showing how computations with appropriate twisted coefficients give fixed point bounds over entire mapping classes.

\begin{thm}
\label{degbound}
The minimum number of fixed points of an area-preserving map $\phi$ in a mapping class $h$ is given by:
$$ \left\{ \begin{array}{cc}
\#\left\{\eta:\ind(\eta)\neq 0\right\}  & h\ \textrm{is periodic or pseudo-Anosov} \\
\#\left\{\eta:\ind(\eta)\neq 0\right\} + A + B & h\ \textrm{is reducible} \end{array} \right.$$
where $A$ is as before and $B$ is the number of fixed annuli.
\end{thm}

Again the upper bound is given by construction. The lower bound comes from ``quantum cup-length'' computations for a certain cohomology operation\footnote{We note that the usual module structure over the quantum cohomology of $\Sigma$ vanishes in the situation of interest.} on the summand of Floer homology corresponding to the given Nielsen class plus a compactness argument of Taubes. This cohomology operation is morally given by counting intersections of holomorphic cylinders with cycles in $H_1(S,\del S)$ for $S$ a subsurface corresponding to the Nielsen class. The fact that this is well defined comes from an understanding of the homotopy type of the component of the space of sections $\Gamma(M_\phi)$ of $M_\phi \rightarrow S^1$ corresponding to the given Nielsen class, plus an algebraic invariance result given in \cite{c} for certain types of cohomology operations on Floer homology. 

The symplectic Floer homology of a symplectomorphism $\phi$ is the homology of a chain complex generated by the fixed points of $\phi$ (see \S\ref{floer} for more details). Thus the rank of Floer homology gives a bound on the number of fixed points of $\phi$, and results on the invariance of Floer homology under deformation of $\phi$ give more general fixed point bounds. A main issue here is that we are interested in fixed point bounds on connected components of $\Symp(X)$ as opposed to on $\Ham(X)$-cosets of $\Symp(X)$. Floer homology is invariant under Hamiltonian perturbations between maps with nondegenerate fixed points, so as long as we can define it on a given $\Ham(X)$-coset, the rank of the Floer homology for any map in the coset is a bound on the number of fixed points for all nondegenerate maps in the coset.

To deal with this, we give a general method for finding fixed point bounds for symplectic mapping classes on monotone symplectic manifolds $(X,\omega)$. We say that $(X,\omega)$ is monotone if $[\omega] \in H^2(X)$ is a positive multiple of $c_1(X)$.\footnote{If the multiple is negative, everything goes through similarly if there are no spheres with chern class above $-n-2$, where $n$ is half the dimension of $X$. Even failing this, we should be able to use virtual moduli space or polyfold methods. These are technical issues, but the statement that $[\omega]$ is some multiple of $c_1(X)$ seems vital to our argument.} The method consists of performing a single Floer homology computation, for a suitable map $\phi$ in the symplectic mapping class, with suitable twisted coefficients. Along the way, we show that Floer homology computations for one suitable map $\phi$ with various twisted coefficients give computations of $HF_*(\psi,\Lambda_\psi)$, the Floer homology of $\psi$ with its natural Novikov coefficients $\Lambda_\psi$, for any $\psi$ in the symplectic mapping class.

We identify a subset of each mapping class which we call \emph{weakly monotone for every Nielsen class} such that Floer homology is defined with any coefficients and is invariant under deformations through such maps. In \cite{c} we have shown that standard form maps on surfaces are in this subset. Let $\omega_\phi$ be the two-form on $M_\phi$ induced by $\omega$ on $\Sigma\times\RR$ and let $c_\phi$ be the first Chern class of the vertical tangent bundle of $M_\phi \rightarrow S^1$.

\begin{defn}
\label{wkmonointro}
A map $\phi: \Sigma \rightarrow \Sigma$ is \emph{weakly monotone}\footnote{We used ``weakly monotone for every Nielsen class'' to refer to the same concept in \cite{c}.} if $[\omega_\phi]$ vanishes on $T_0(M_\phi)$, where $T_0(M_\phi) \subset H_2(M_\phi;\RR)$ is generated by tori $T$ with $c_\phi(T) = 0$ such that $\pi|_T: T \rightarrow S^1$ is a fibration with fiber $S^1$, where the map $\pi: M_\phi \rightarrow S^1$ is the projection.
\end{defn}

We define a Flux map $\Flux(\phi): N_h\rightarrow \RR$. Here $N_h$ denotes the image under the map $H_2(M_\phi) \rightarrow \ker(\phi_*-id) \subset H_1(\Sigma)$ of $\ker(c_\phi) \subset H_2(M_\phi)$. Let $N_h'$ be the image under the map $H_2(M_\phi) \rightarrow \ker(\phi_*-id) \subset H_1(\Sigma)$ of the subset $T_0(M_\phi) \subset \ker(c_\phi) \subset H_2(M_\phi)$.

\begin{thm}
\label{novikovthm}
Let $X$ be a monotone symplectic manifold, $h$ a symplectic mapping class on $X$, and $\phi \in h$.
\be
\ii $\Flux(\phi): N_h/tors \rightarrow \RR$ is well-defined.
\ii $\Flux(\phi)|_{N_h'/tors} = 0$ if and only if $\phi$ is weakly monotone for every Nielsen class.
\ii Let $\phi$ be such that $\Flux(\phi)|_{N_h} = 0$. Then the rank of $HF_*(\phi; Q(\ZZ/2[N]))$ gives a lower bound on the rank $HF_*(\psi; \Lambda_\psi)$ for $\psi$ any map in the class $h$, where $\Lambda_\psi$ is the Novikov ring over which $HF_*(\psi)$ is naturally defined, and thus gives a lower bound on the number of fixed points of a map in the mapping class with nondegenerate fixed points.
\ee
\end{thm}

Here $Q(\ZZ/2[N_h'/tors])$ is the quotient field of the group ring of $N_h'/tors$ over $\ZZ/2$. The Novikov ring $\Lambda_\psi$ is a ring over which $HF_*(\psi)$ can be defined, with generators given by fixed points, even if $\psi$ is not weakly monotone. The main idea is that in some cases the field $Q(\ZZ/2[N_h])$ injects into $\Lambda_\psi$, in which case we have a field extension and the ranks of the homology are the same. When we don't have an injection, we can extend $\Lambda_\psi$ to a larger Novikov ring into which $Q(\ZZ/2[N_h])$ does inject. When homology is computed over the larger Novikov ring, the rank can only decrease, giving the lower bound in the theorem. Yi-Jen Lee's bifurcation analysis \cite{l1,l2} is vital to this argument, giving a way to compare the Floer homology of maps with proportional fluxes. We note that Lee and Taubes \cite[Corollary 6.6]{lt} have a similar theorem for periodic Floer homology coming from their isomorphism with Seiberg-Witten Floer homology. Lee also informed us that she was independently aware of the existence of results such as Theorem \ref{novikovthm}.

\subsection{Organization of the paper}

In \S\ref{floer} we review Floer homology and Nielsen classes; develop general twisted coefficients for Floer homology; review Novikov rings; and discuss invariance results due to Yi-Jen Lee \cite{l1,l2}.

In \S\ref{bounds} we give a general method for finding fixed point bounds on a monotone symplectic manifold using a Floer homology computation for any weakly monotone symplectomorphism with a particular choice of twisted coefficients.

In \S\ref{surfacebounds} we carry out this method in the case of surface symplectomorphisms to give a lower bound, using computations from \cite{c}. We obtain an equal upper bound by explicit constructions.

In \S\ref{degenerate}, we give fixed point bounds for surface symplectomorphisms with possibly degenerate fixed points. We use a certain cohomology operation in place of the quantum cap product, which vanishes in the situation of interest, to give cup-length-type bounds. An equal upper bound is again given by explicit constructions.

\subsection{Acknowledgments}

This paper is adapted from part of my thesis \cite{c2} at UC Berkeley under Michael Hutchings, who gave invaluable advice and support. I'd also like to thank MIT for their hospitality during my last year of graduate school, and Denis Auroux, for his advising and helpful discussions.

\section{Floer homology, Nielsen classes, twisted coefficients, Novikov rings, and bifurcation analysis}
\label{floer}

\subsection{Review of Floer theory and monotonicity}

We provide a brief summary. For a more complete discussion, see \cite{ds, s2, c, per1}.

Let $\Sigma$ be a compact, connected, oriented surface, possibly with boundary, of negative Euler characteristic. Let $\omega$ be a symplectic form (i.e. an area form) on $\Sigma$. Let $\phi \in \Symp^\del(\Sigma,\omega)$, the space of symplectomorphisms (i.e. area-preserving diffeomorphisms) with no fixed points on the boundary. We consider the mapping torus of $\phi$: $$M_\phi = \frac{\RR \times \Sigma}{(t+1,x)\simeq(t,\phi(x))}.$$ Note that this is a $\Sigma$-bundle over $S^1$ with projection $\pi: M_\phi \rightarrow \RR/\ZZ = S^1$.

Let $\Gamma(M_\phi)$ denote the space of smooth sections of $\pi: M_\phi \rightarrow S^1$. Note that a fixed point $x\in \Sigma$ of $\phi$ can be interpreted as a constant section $\gamma_x$. Let $\mathcal{J}$ denote the space of almost complex structures $J$ on $\RR \times M_\phi$ which are $\RR$-invariant, preserve the vertical tangent bundle of $\pi: \RR \times M_\phi \rightarrow \RR \times S^1$, and for which $\pi$ is $J-j$-holomorphic, given the standard complex structure $j$ on the cylinder $\RR\times S^1$.

Suppose $\phi$ has nondegenerate fixed points, in the sense that $d\phi$ does not have $1$ as an eigenvalue at any fixed point. Let $P_{x,y}\Gamma(M_\phi)$ denote the space of paths from $\gamma_x$ to $\gamma_y$ in $\Gamma(M_\phi)$. Let $C\in \pi_0(P_{x,y}\Gamma(M_\phi))$. Given a generic $J \in \mathcal{J}$, the moduli space $\mathcal{M}(\phi,x,y,C)$ of holomorphic sections $\RR\times S^1 \rightarrow \RR\times M_\phi$ in the homotopy class $C$ is smooth and compact, of dimension $\ind(C)$ (see \cite{fhs}, \cite{gr}).

For a statement of the index formula, see \cite{rs}, or \cite[\S2.2]{c} for a discussion tailored to the current setting. Other than computations that we will cite below from \cite{c}, we have need only of the change of homology class formula. In preparation, we define the cohomology class $c_\phi \in H^2(M_\phi) = H^2(\RR\times M_\phi)$ to be the first chern class of the vertical tangent bundle to the projection $\pi$.

\begin{prop}[{\cite{rs}}]
Let $C, C' \in \pi_0(P_{x,y}\Gamma(M_\phi))$. Then $\ind(C) - \ind(C') = 2\left<c_\phi,[C-C']\right>$.
\end{prop}

Let $\omega_\phi$ denote the cohomology class in $H^2(M_\phi) = H^2(\RR\times M_\phi)$ of the vertical area form on $M_\phi$.

\begin{prop}[{Gromov Compactness, \cite{gr}}]
Let $I \subset \pi_0(P_{x,y}\Gamma(M_\phi))$. The union $\cup_{C\in I} \mathcal{M}(\phi,x,y,C)$ is compact if the set of values $$\left\{\omega_\phi(C-C')| C,C' \in \pi_0(P_{x,y}\Gamma(M_\phi)\right\}$$ is bounded.
\end{prop}

We have the following conditions:

\begin{defn}
\label{monofloer}
A map $\phi \in \Symp^\del(\Sigma,\omega)$ is \emph{monotone} if $\omega_\phi$ vanishes on the kernel of $c_\phi$.
\end{defn}

\begin{defn}
\label{wkmonofloer}
A map $\phi \in \Symp^\del(\Sigma,\omega)$ is \emph{weakly monotone} if $\omega_\phi$ vanishes on $T_0(M_\phi)$, where $T_0(M_\phi) \subset H_2(M_\phi;\RR)$ is generated by tori $T$ with $c_\phi(T) = 0$ such that $\pi|_T: T \rightarrow S^1$ is a fibration with fiber $S^1$, where the map $\pi: M_\phi \rightarrow S^1$ is the projection.
\end{defn}

Under either of these conditions, the moduli space $$\mathcal{M}_1(\phi,x,y) = \cup_{C : \ind(C) = 1} \mathcal{M}(\phi,x,y,C)$$ is compact. Thus we may define they symplectic Floer homology $HF_*(\phi)$ of $\phi$ with coefficients in $\ZZ/2$ to be the homology of the $\ZZ/2$-graded chain complex $CF(\phi) = \oplus_{x \in \mathrm{Fix}(\phi)} \ZZ/2 \cdot x$ with differential $$\del x = \sum_y \# \left(\mathcal{M}(\phi,x,y)/\RR\right) \cdot y,$$ where the $\RR$ action is by translation in the $\RR$-direction in $\RR\times M_\phi$, and the $\ZZ/2$-grading of a fixed point $x$ is given by the sign of $\det(1-d\phi)$. The homology of this chain complex is invariant under deformations through monotone or weakly monotone maps (see \cite{s2, c}). We recall:

\begin{prop}[{\cite{s2} for monotone, \cite{c} for weakly monotone}]
The space of monotone (resp. weakly monotone for every Nielsen class) maps $\phi \in \Symp^\del(\Sigma,\omega)$ is homotopy equivalent to $\Diff^\del(\Sigma)$ under the inclusion map. In particular, the space of those in each mapping class is connected.
\end{prop}

\subsection{Nielsen classes and Reidemeister trace}

There is a topological separation of fixed points due to Nielsen: given a fixed point $x$ of $\phi$, we obtain an element $[\gamma_x] \in \pi_0(\Gamma(M_\phi)$, the homotopy class of $\gamma_x$ in the space of sections $\Gamma(M_\phi)$. The chain complex $CF_*(\phi)$ defined above, as well as all the variants to be defined below, split into direct summands for each Nielsen class $\eta \in \pi_0(\Gamma(M_\phi)$: $$HF_*(\phi) = \sum_{\eta\in\pi_0(\Gamma(M_\phi)} HF_*(\phi,\eta).$$

We note that in the case of surfaces of negative Euler characteristic, Nielsen classes are well-defined on entire mapping classes (i.e. there is no monodromy). Thus the rank of $HF_*(\phi,\eta)$ is a lower bound on the number of fixed points in the Nielsen class $\eta$ on the e.g. weakly monotone subset of the mapping class.

A simpler lower bound on the number of nondegenerate fixed points in a Nielsen class $\eta \in \pi_0(\Gamma(M_\phi)$ is given by the absolute value of the index of the Nielsen class $$\ind(\eta) = \sum_{x: [\gamma_x]=\eta} \ind(x).$$ We note that $\ind(\eta)$ is the Euler characteristic of $HF_*(\phi,\eta)$. A lower bound for the number of fixed points of a map $\psi$ with nondegenerate fixed points in the same mapping class as $\phi$ is thus given by the \emph{Reidemeister trace} $$R(\phi) = \sum_{\eta \in \pi_0(\Gamma(M_\phi)} |\ind(\eta)|.$$

For maps with possibly degenerate fixed points, the bound from Nielsen theory is the number of \emph{essential} Nielsen classes. That is, the number of Nielsen classes $\eta$ for which $\ind(\eta) \neq 0$.

\subsection{Twisted coefficients}
\label{twisted}

We have a groupoid $\GG_\phi$ whose objects are fixed points and whose morphisms from $x$ to $y$ are given by homotopy classes of paths from $x$ to $y$ in $\Gamma(M_\phi)$. Given a ring $R$, we can take the groupoid algebroid $R(\GG_\phi)$, a category enriched over $R$-modules with the same objects as $\GG_\phi$ and with morphisms given by free $R$-modules generated by morphisms in $\GG_\phi$. Note that we have a homomorphism $\ind: \GG_\phi \rightarrow \ZZ$ given by the index.

If we have a representation $\rho$ from $R(\GG_\phi)$ to an $R$-module $M$, we can define \emph{Floer homology with coefficients in $\rho$ (or $M$ if $\rho$ is understood))} as the homology of a chain complex over $M$ with generators the fixed points of $\phi$ and differential given by $$\del x = \sum_y \sum_{C\in \MM_1(x,y)} \rho(C)\cdot y,$$ where $[C]$ is the homotopy class of the path in $\Gamma(M_\phi)$ associated to the flow line $C$. This is defined for weakly monotone $\phi$ for arbitrary $\rho$ (because then $\MM_1(\phi,x,y)$ is compact). In \S\ref{novikov} we discuss representations $\Lambda_\phi$ suitable for each $\phi$.

We will typically suppress the ring $R$, which throughout this paper may be assumed to be $\ZZ/2$.

The \emph{standard representation} $\rho_{st}$ into the group ring of $H_1(\Gamma(M_\phi))$ is defined as follows. For every pair of fixed points $x, y$ and every index $i$, we choose a path $C^i_{x,y}$ in $\Gamma(M_\phi)$ between them of index $i$ if possible. We require that $C^i_{x,y}\cdot C^j_{y,z} \simeq C^{i+j}_{x,y}$, that $C^i_{x,y} = -C^{-i}_{y,x}$, and that $C^0_{x,x} = *$. Here $\cdot$ signifies appending paths, $-$ signifies reversal of a path, and $*$ signifies the constant path. Then $\rho_{st}([C])$, for $C$ a path from $x$ to $y$, is defined to be $[C\cdot C^{-\ind(C)}_{y,x}] \in H_1(\Gamma(M_\phi))$. Note that in fact this lies in $\ker(c_\phi)\subset H_1(\Gamma(M_\phi))$. We have made choices, but the resulting Floer homology is well-defined up to a change of basis.

We typically compose this with the map $H_1(\Gamma(M_\phi)) \rightarrow H_2(M_\phi)$ to get what we call \emph{fully twisted coefficients}, over $\ZZ/2[H_2(M_\phi)]$. If we desire to have field coefficients, we may for example take the quotient field of the group ring of $H_2(M_\phi)/tors$.

\subsection{Novikov rings}
\label{novikov}

\begin{defn}
\label{novikovdef}
For $R$ a ring, $G$ an abelian group and a homomorphism $N: G \rightarrow \RR$, the \emph{Novikov ring} $\Nov_R(G,N)$ is defined to be the ring whose elements are formal sums $\sum_{g\in G} a_g \cdot g$, where $a_g \in R$ are such that, for each $r \in \RR$, only finitely many of the $a_g$ with $N(g) < r$ are nonzero.
\end{defn}

There is also the universal Novikov ring:

\begin{defn}
\label{universal}
For $R$ a ring, the \emph{universal Novikov ring} $\Lambda_R$ is defined to be the ring whose elements are formal sums $\sum_{r\in R} a_r \cdot T^r$, where $a_r \in R$ are such that, for each $s \in \RR$, only finitely many of the $a_r$ with $r < s$ are nonzero.
\end{defn}

Note that this is a field if $R$ is a field. We have maps $\Nov_R(G,N) \rightarrow \Lambda_R$ given by $g \mapsto T^{N(g)}$. If $R$ is a field and $N$ is injective, this is an extension of fields.

We define a representation, the \emph{natural Novikov coefficients for $\phi$}, denoted $\Lambda_\phi$, as follows. We have a representation into $\Nov(\ker(c_\phi), \omega_\phi, \ZZ/2)$, where $c_\phi, \omega_\phi : H_1(\Gamma(M_\phi)) \rightarrow \RR$ are defined as in \S\ref{floer}. The representation is defined in the same manner as the standard representation; we have simply taken a submodule $\ker(c_\phi) \subset H_1(\Gamma(M_\phi)$ in which the image must lie and allowed certain infinite sums. We further compose this with the map to $\Lambda_{\ZZ/2}$. This all is simply to say we take $$\rho(C) = T^{\omega_\phi(C\cdot C^{-ind(C)}_{y,x})}.$$ We may define $HF_*(\phi,\Lambda_\phi)$, the Floer homology of $\phi$ with coefficients in $\Lambda_\phi$ for any $\phi$. The point is simply that, while we may not have finiteness for $\MM_1(x,y)$, we do have that there are only finitely many $C\in \MM_1(x,y)$ with $\omega_\phi(C) < r$ for any given $r\in \RR$, by Gromov compactness, and this is precisely what is required.

\subsection{Bifurcation analysis of Yi-Jen Lee}
\label{bif}

Yi-Jen Lee \cite{l1,l2} has worked out a general bifurcation argument for what she calls Floer-type theories.
Michael Usher \cite{u} has a nice summary of the invariance result this gives (which Lee conjectured in an earlier paper \cite[Eqn 3.2]{l} but did not explicitly state as a Theorem in \cite{l1,l2}) and its algebraic aspects. We have adapted the statement to our setting involving representations $\rho$ of $\ZZ/2(\GG_\phi)$.

\begin{thm}[{\cite[Theorem 3.6, due to Lee (\cite{l1},\cite{l2})]{u}}]
\label{invariance}
Suppose $(\Sigma, \omega)$ is a symplectic manifold with $\pi_2(\Sigma) = 0$. Let $\phi_r: \Sigma \rightarrow \Sigma$ be a smooth family of symplectomorphisms and $J_r = \{J_t\}_r$ a smooth family (of $1$-periodic families) of almost complex structures such that $(\phi_0,J_0)$ and $(\phi_1,J_1)$ are generic. Let $N = \ker(c_\phi): H_1(\Gamma(M_{\phi_r})_\eta) \rightarrow \RR$ (this is independent of $r$). 
\be
\ii Suppose $\omega_{\phi_r}|_N = f(r) \omega_{\phi_0}|_N$, for $f(r) \in \RR$.
Then $HF_*(\phi_0,\eta, J_0; \Lambda_{\phi_0}) \cong HF_*(\phi_1,\eta, J_1; \Lambda_{\phi_0})$.\qed
\ii Suppose $\omega_{\phi_r}|_N = 0$ for all $r$. Then $HF_*(\phi_0,\eta, J_0) \cong HF_*(\phi_1,\eta, J_1)$. In fact,  $HF_*(\phi_0,\eta, J_0;\rho) \cong HF_*(\phi_1,\eta, J_1;\rho)$ for any representation $\rho$.
\ee
\end{thm}

We point out that in item (1) we are indeed taking the Floer homology of $\phi_1$ with coefficients given by the representation $\Lambda_{\phi_0}$. What this means is that we use $\omega_{\phi_0}$ to determine the power of $T$ in $\rho(C)$.

\section{Bounds on fixed points in symplectic mapping classes}
\label{bounds}

\subsection{Representations and flux}
The goal of this section is to give, for $(X,\omega)$ monotone, a lower bound on fixed points for symplectomorphisms in a symplectic mapping class $h \in \Gamma_\Symp(X,\omega)$ in terms of the rank of $HF_*(\phi,\rho)$ for one suitable choice of a pair $(\phi,\rho)$ with $\phi \in h$ and $\rho$ a representation of $\ZZ/2(\GG_\phi)$. In this subsection we define the representations $\rho_m$ for which it will be shown that $(\phi,\rho_m)$ is such a pair for any monotone $\phi$ and $\rho_{wm}$ for which it will be shown that $(\phi,\rho_{wm})$ is such a pair for any weakly monotone $\phi$.

\subsection*{Monotone case}

For $(X,\omega)$ monotone, we have a representation $\rho_m$ that works for any monotone $\phi$. The representation $\rho_m$ is defined as follows:

We have the representation $\rho_{st}$ into $\ZZ/2[H_1(\Gamma(M_\phi))]$. We first compose with $H_1(\Gamma(M_\phi)) \rightarrow H_2(M_\phi)$, noting that the image lies in the kernel of $c_\phi: H_2(M_\phi) \rightarrow \RR$. Now we compose with the map to $H_1(X)$ in the long exact sequence for the mapping torus, a part of which is $$H_2(X) \stackrel{i}{\rightarrow} H_2(M_\phi) \stackrel{\del}{\rightarrow} \ker(\phi_*-\id) \subset H_1(X).$$ We denote the image of $\ker(c_\phi)$ in $H_1(X)$ by $N_h$ (this depends only on the mapping class $h$). Finally, we mod out by torsion and then take the quotient field, so that our coefficients lie in the field $$Q(\ZZ/2[N_h/tors]),$$ where $Q$ denotes taking the quotient field. Essentially we have taken fully twisted field coefficients, but we've been careful about where the image lies so that we can define flux in this context, which both will be useful in specifying which maps $\phi$ are suitable to work with and will be a useful tool in the proof that is to come.

\begin{lemmadef}
\label{flux}
The \emph{flux} of a symplectomorphism $\phi$, denoted $$\Flux(\phi): N_h/tors \rightarrow \RR,$$ defined as $$\Flux(\phi)(\gamma) = \omega_\phi(C_\gamma),$$ where $C_\gamma \in H_2(M_\phi)$ is such that $\del(C_\gamma) = \gamma$ and $c_\phi(C_\gamma) = 0$, is well-defined.
\end{lemmadef}

\proof
First, we note that such a $C_\gamma$ exists because $N_h$ is the image under $\del$ (in the long exact sequence) of $\ker(c_\phi)$. Next we note that $\omega_\phi(C_\gamma)$ is well-defined because if we take any other $C'_\gamma$ such that $\del(C'_\gamma) = \gamma$ and $c_\phi(C'_\gamma) = 0$, the difference $C_\gamma - C'_\gamma$ is in the kernel of $c_\phi$ and, being in the kernel of $\del$, must also come from an element of $H_2(\Sigma)$, which we denote $B$. We have that $c_\phi \circ i : H_2(X) \rightarrow \ZZ$ agrees with $c_1(X) : H_2(X) \rightarrow \ZZ$. Thus $c_1(X)(B) = 0$. By monotonicity of $X$, this implies that $\omega(B) = 0$.  Because $\omega_\phi \circ i: H_2(X) \rightarrow \RR$ agrees with $\omega: H_2(X) \rightarrow \RR$ as well, we conclude that $\omega_\phi(C_\gamma) - \omega_\phi(C'_\gamma) = 0$. Finally, we remark that this is a homomorphism, and thus any torsion in $N_h$ must map to $0 \in \RR$.
\qed

\begin{lemma}
For a symplectomorphism $\phi$, we have: $$\phi\ \textrm{is monotone}\ \iff \Flux(\phi) = 0.$$
\end{lemma}

\proof
The statement $\Flux(\phi) = 0$ is equivalent to $\omega_\phi$ vanishing on the kernel of $c_\phi$. This implies that $\omega_\phi$ is some multiple of $c_\phi$.
\qed

\subsection*{Weakly monotone case}

With $(X,\omega)$ monotone, we also give a representation $\rho_{wm}$ that works for any weakly monotone $\phi$. The representation $\rho_{wm}$ is defined as follows:

Let the image of the map $H_1(\Gamma(M_\phi)) \rightarrow H_2(M_\phi)$ be denoted $T(M_\phi)$. This is generated by $2$-tori in $M_\phi$ standardly fibering over $S^1$. When composing the standard representation into $\ZZ/2[H_1(\Gamma(M_\phi)]$ with the map $H_1(\Gamma(M_\phi)) \rightarrow H_2(M_\phi)$, the image thus lies not only in $\ker(c_\phi)$ but in $\ker(c_\phi|_{T(M_\phi)})$. In fact it lies moreover in $T_0(M_\phi)$, generated by $2$-tori $T$ in $M_\phi$ standardly fibering over $S^1$ such that $c_\phi(T) = 0$. We denote the image of $T_0(M_\phi)$ under $\del$ in the long exact sequence as $N_h' \subset N_h \subset \ker(\phi_*-\id) \subset H_1(X)$. As before, we mod out by torsion and then take the quotient field, so our coefficients lie in the field $$Q(\ZZ/2[N_h'/tors]).$$

We restrict $Flux(\phi)$ to $N_h'/tors \subset N_h/tors$.

\begin{lemma}
For a symplectomorphism $\phi$, we have: $$\phi\ \textrm{is weakly monotone for every Nielsen class $\eta$}\ \iff \Flux(\phi)|_{N_h'/tors} = 0.$$
\end{lemma}

\proof
The statement $\Flux(\phi)_{N_h'/tors} = 0$ is equivalent to $\omega_\phi$ vanishing on $T_0(M_\phi)$.  Let $T_0(M_\phi)_\eta$ be generated by $2$-tori $T$ in $M_\phi$ standardly fibering over $S^1$ such that $c_\phi(T) = 0$ and which moreover have a section in Nielsen class $\eta$. Then the vanishing of $\omega_\phi$ on $T_0(M_\phi)$ is equivalent to $\omega_\phi$ vanishing on $T_0(M_\phi)_\eta$ for every $\eta$, because every such torus has a section, which lies in some Nielsen class (i.e. component of $\Gamma(M_\phi)$).
\qed

\begin{cor}
The space of maps in the mapping class $h$ which are weakly monotone for every Nielsen class $\eta$, $\Symp_h^{wm}(X,\omega)$, is homotopy equivalent under the inclusion to $\Symp_h(X,\omega)$.\qed
\end{cor}

\subsection{Fixed point bounds}

\begin{thm}
\label{generalbounds}
Let $X$ be a monotone symplectic manifold, $h$ a symplectic mapping class on $X$, and $\phi \in h$. Let $\phi$ be such that $\Flux(\phi)|_{N_h'} = 0$. Then the rank of $$HF_*(\phi; Q(\ZZ/2[N_h'/tors]))$$ gives a lower bound on the rank $HF_*(\psi; \Lambda_\psi)$ for $\psi$ any map in the class $h$.
\end{thm}

\begin{cor}
\label{boundcor}
Let $X$ be a monotone symplectic manifold, $h$ a symplectic mapping class on $X$, and $\phi \in h$. Let $\phi$ be such that $\Flux(\phi)|_{N_h'} = 0$. Then the rank of $$HF_*(\phi; Q(\ZZ/2[N_h'/tors]))$$ gives a lower bound on the number of fixed points of any map in the mapping class with nondegenerate fixed points.\qed
\end{cor}

\begin{remark}
The corresponding versions with $\phi$ such that $\Flux(\phi) = 0$ and coefficients in $Q(\ZZ/2[N_h/tors])$ also hold, either with the same proofs or as a formal consequence.
\end{remark}

Theorem \ref{generalbounds} follows from the following three Lemmas.

\begin{lemma}
\label{lemma1}
Let $\phi \in h$ be such that $\Flux(\phi)|_{N_h'} = 0$ and let $\psi \in h$. Then $HF(\psi,\Lambda_\psi) \cong HF(\phi,\Lambda_\psi)$.
\end{lemma}

\proof
We use (1) from Theorem \ref{invariance}. Consider a smooth family $\psi_r$ of symplectomorphisms from $\psi_0 = \psi$ to $\psi_1 = \phi$ with $\Flux(\psi_r)|_{N_h'} = (1-r)\Flux(\psi)|_{N_h'}$. Then (for a generic family of almost complex structures) the conditions of (1) in Theorem \ref{invariance} are met with $f(r) = 1-r$. Summing over Nielsen classes gives the result.
\qed

\begin{lemma}
\label{lemma2}
With $\phi$, $\psi$ as above, the rank of $HF(\phi,\Lambda_\psi)$ equals the rank of $$HF(\phi,\ZZ/2[N_h'/tors/\ker(\Flux(\psi)|_{N_h'})]).$$
\end{lemma}

\proof
By the rank over $\ZZ/2[N_h'/tors/\ker(\Flux(\psi)|_{N_h'})]$, we mean the rank over $Q(\ZZ/2[N_h'/tors/\ker(\Flux(\psi)|_{N_h'})])$, the quotient field. We have field extensions \beqa Q(\ZZ/2[N_h'/tors/\ker(\Flux(\psi)|_{N_h'})]) &\hookrightarrow& \\ \Nov\left(N_h'/tors/\ker(\Flux(\psi)|_{N_h'}),\Flux(\psi)|_{N_h'}\right) &\hookrightarrow& \Lambda_\psi.\eeqa The first is by allowing some infinite sums, and the second is the map discussed in \S\ref{novikov}. Field extensions are flat, so the ranks are equal.
\qed

\begin{lemma}[{(\cite[Appendix C, with T.-Q.-T. L\^e]{lo})}]
\label{lemma3}
Let $k$ be a field and let $C_*$ be a chain complex over $k_n := k[t_1,\ldots,t_n]$. Consider $k_n$ as an augmented $k_m$-algebra (for some $0 < m < n$) with augmentation sending $t_i$ to $t_i$ for $i \leq m$ and $t_i$ to $1$ for $i > m$. Then 
$$\textrm{rank of}\ H_*(C_*\otimes_{k_n} k_m) \leq \textrm{rank of}\ H_*(C_*).$$\qed
\end{lemma}

\textit{Proof of Theorem \ref{generalbounds}:}\ \ \ \ 
Lemma \ref{lemma1} implies that in particular the rank of $HF(\psi,\Lambda_\psi)$ equals the rank of $HF(\phi,\Lambda_\psi)$. By Lemma \ref{lemma2} this rank is equal to the rank of $HF(\phi,\ZZ/2[N_h'/tors/\ker(\Flux(\psi)|_{N_h'})])$. We note that
\beqa &CF(\phi,\ZZ/2[N_h'/tors/\ker(\Flux(\psi)|_{N_h'})]) \simeq& \\
&CF(\phi,\ZZ/2[N_h'/tors]) \otimes_{\ZZ/2[N_h'/tors]} \ZZ/2[N_h'/tors/\ker(\Flux(\psi)|_{N_h'})]&\eeqa
with the adjunction map sending elements of the kernel of $\Flux(\psi)|_{N_h'}$ to one. The spaces $N_h'/tors$ and $N_h'/tors/\ker(\Flux(\psi)|_{N_h'})$ are vector spaces over $\RR$ and as such their group rings over $\ZZ/2$ and the aforementioned adjunction map are as in Lemma \ref{lemma3}, where $n$ and $m$ are their dimensions repsectively, with $k = \ZZ/2$. Thus by Lemma \ref{lemma3} the rank of the homology of the displayed complexes is less than or equal to the rank of $HF(\phi,\ZZ/2[N_h'/tors])$, which is the same as the rank of $HF(\phi,Q(\ZZ/2[N_h'/tors]))$, giving the result.
\qed

\section{Bounds on fixed points for surface symplectomorphisms}
\label{surfacebounds}

\subsection{Standard form maps and Nielsen classes}
\label{red2}

In this section we describe standard form maps and their Nielsen classes. For the identity mapping class, a standard form map is a small perturbation of the identity map by the Hamiltonian flow associated to a Morse function for which the boundary components are locally minima and maxima. Every fixed point is in the same Nielsen class. This Nielsen class has index given by the Euler characteristic of the surface. For non-identity periodic mapping classes, a standard form map is an isometry with respect to a hyperbolic structure on the surface with geodesic boundary. Every fixed point is in a separate Nielsen class and each of the Nielsen classes for which there is a fixed point has index one. For a pseudo-Anosov mapping classes, a standard form map is a symplectic smoothing of the singularities and boundary components of the standard singular representative. Each singularity has a number $p\geq 3$ of prongs and each boundary component has a number $p \geq 1$ of prongs. If a singularity or boundary component is (setwise) fixed, it has some fractional rotation number modulo $p$. There is a separate Nielsen class for every smooth fixed point, which is of index one or minus one; for every fixed singularity, which when symplectically smoothed gives $p-1$ fixed points all of index minus one if the rotation number is zero modulo $p$ or one fixed point of index one otherwise \cite{c}; and for every fixed boundary component with rotation number zero modulo $p$, which when symplectically smoothed gives $p$ fixed points all of index minus one.

From this discussion, we see that for non-identity periodic and pseudo-Anosov mapping classes, the standard form map is such that all fixed points are nondegenerate and, for every Nielsen class $\eta$, the number of fixed points in $\eta$ is $|\ind(\eta)|$. We now turn to reducible maps and the identity map.

By Thurston's classification (see \cite{t} and \cite{flp}; also cf. \cite[Definition 8]{g} and \cite[Definition 4.6]{c}), in a reducible mapping class $h$, there is a (not necessarily smooth)  map $\phi$ which satisfies the following:

\begin{defn}
\label{reduciblestandard}
A reducible map $\phi$ is in \emph{standard form} if there is a $\phi$-and-$\phi^{-1}$-invariant finite union of disjoint noncontractible (closed) annuli $N \subset \Sigma$ such that:
\be
\ii \label{twist} For $A$ a component of $N$and $\ell$ the smallest positive integer such that $\phi^\ell$ maps $A$ to itself, the map $\phi^\ell|_A$ is either a \emph{twist map} or a \emph{flip-twist map}. That is, with respect to coordinates $(q,p)\in[0,1]\times S^1$, we have one of \beqa(q,p)&\mapsto& (q,p-f(q)) \qquad \textrm{(twist map)} \\ (q,p)&\mapsto& (1-q,-p+f(q)) \qquad \textrm{(flip-twist map)},\eeqa
where $f : [0,1]\rightarrow \RR$ is a strictly monotonic smooth map. We call the (flip-)twist map \emph{positive} or \emph{negative} if $f$ is increasing or decreasing, respectively. Note that these maps are area-preserving.
\ii Let $A$ and $\ell$ be as in (\ref{twist}). If $\ell = 1$ and $\phi|_N$ is a twist map, then $\Im(f)\subset[0,1]$. That is, $\phi|_{\textrm{int}(A)}$ has no fixed points. (If we want to twist multiple times, we separate the twisting region into parallel annuli separated by regions on which the map is the identity.) We further require that parallel twisting regions twist in the same direction.
\ii For $S$ a component of $\Sigma \backslash N$ and $\ell$ the smallest integer such that $\phi^\ell$ maps $A$ to itself, the map $\phi^\ell|_S$ is area-preserving and is either isotopic to the \emph{identity},  \emph{periodic}, or \emph{pseudo-Anosov}. In these cases, we require the map to be in standard form as above.
\ee
\end{defn}

We review Nielsen classes of fixed points of standard form reducible maps, as discussed in \cite{c}. The fixed points of our standard form reducible maps are as follows:

\bi
\ii (Type Ia) The entire component of fixed components $S$ of $\Sigma\backslash N$ with $\chi(S) < 0$.
\ii (Type Ib) The entire component of fixed components $S$ of $\Sigma\backslash N$ with $\chi(S) = 0$. These are annuli and only occur when we have multiple parallel Dehn twists.
\ii (Type IIa) Fixed points $x$ of periodic components $S$ of $\Sigma\backslash N$ with $\chi(S) < 0$ which are setwise fixed by $\phi$. These are each index one.
\ii (Type IIb) Fixed points $x$ of flip-twist regions. These are each index one. Note that each flip-twist region has two fixed points.
\ii (Type III) Fixed points $x$ of pseudo-Anosov components $S$ of $\Sigma\backslash N$ which are setwise fixed by $\phi$. These come in 4 types (note that there are no fixed points associated to a rotated puncture):
\bi
\ii (Type IIIa) Fixed points which are not associated with any singularity or puncture (i.e. boundary component) of the pre-smoothed map. These are index one or negative one.
\ii (Type IIIb-$p$) Fixed points which come from an unrotated singular point with $p$ prongs. There are $p-1$ of these for each such, each index negative one.
\ii (Type IIIc) Fixed points which come from a rotated singular point. These are each index one.
\ii (Type IIId-$p$) Fixed points which come from an unrotated boundary component with $p$ prongs. There are $p$ for each such, each index negative one.
\ei
\ei

In \cite{c}, adapting the work of \cite{jg} to the area-preserving case, we show that we have a separate Nielsen class for every component of Type Ia or Type Ib, for every single fixed point of Type IIa, IIb, IIIa, or IIIc, and for every unrotated singular point of the pre-smoothed map for Type IIIb (i.e. the collection of fixed points associated to a single unrotated singular point are all in the same Nielsen class).

Type IIId fixed points associated to the same boundary component are in the same Nielsen class. They may also be in the same Nielsen class as fixed points of the region they abut if that region is also of type IIId or is of type Ia (they cannot abut regions of type Ib). In the former case this Nielsen class is again separate from all others already specified, and has combined index $-p-q$ from the $p+q$ index negative one fixed points. In the latter case, we have already stated that this Nielsen class is separate from all others already specified. Thus the combined index of the Nielsen class $\eta$ associated to a fixed component of $S$ is $$\ind(\eta) = \chi(S) - \sum_{C \in\pi_0(\del S):\ \textrm{the other side of}\ C\ \textrm{is of type IIId-}p} p.$$
Finally, the index of a fixed component of type Ib is zero.

\subsection{Floer homology with twisted coefficients}

We compute the Floer homology $HF_*(\phi,Q(\ZZ/2[N_h'/tors]))$ for $\phi$ a standard form reducible map. This splits into a direct sum over Nielsen classes. For Nielsen classes not associated to a fixed component, $\phi$ has fixed points all of the same index and thus there are no flow-lines at all and the differential vanishes.

For fixed components $S$, possibly abutting type IIId fixed boundary components of pseudo-Anosovs, we perturb with the Hamiltonian flow of a small Morse-Smale function that patches together with the perturbation on a neighborhood of any fixed boundary components, which is given by a Hamiltonian flow of a Morse-Smale function with $p$ saddle points; see \cite{c} for details. In \cite{c} we show that the flow-lines we get between fixed points in the Nielsen class $\eta$ corresponding to this fixed component are only those which correspond to Morse flow-lines.

We are interested in the rank of the summand of $HF_*(\phi,Q(\ZZ/2[N_h'/tors]))$ corresponding to such a Nielsen class $\eta$. The key is to understand the extrema of the Morse-Smale function. If the component has boundaries which rotate in opposite directions, we may choose the Morse-Smale function to have no extrema, and then there are $|\ind(\eta)|$ fixed points in the Nielsen class all of the same index and thus no flow-lines. If the component has boundaries rotating all in the same direction, we may choose a Morse-Smale function with one extremum. There are $|\ind(\eta)|+2$ fixed points in the Nielsen class $\eta$. Finally, if there is no boundary, we may choose a Morse-Smale function with two extrema. There are $|\ind(\eta)|+4$ fixed points in the Nielsen class $\eta$. In these latter two cases, we must further understand the flow-lines.

\begin{figure}
\begin{center}
\includegraphics[scale=.25]{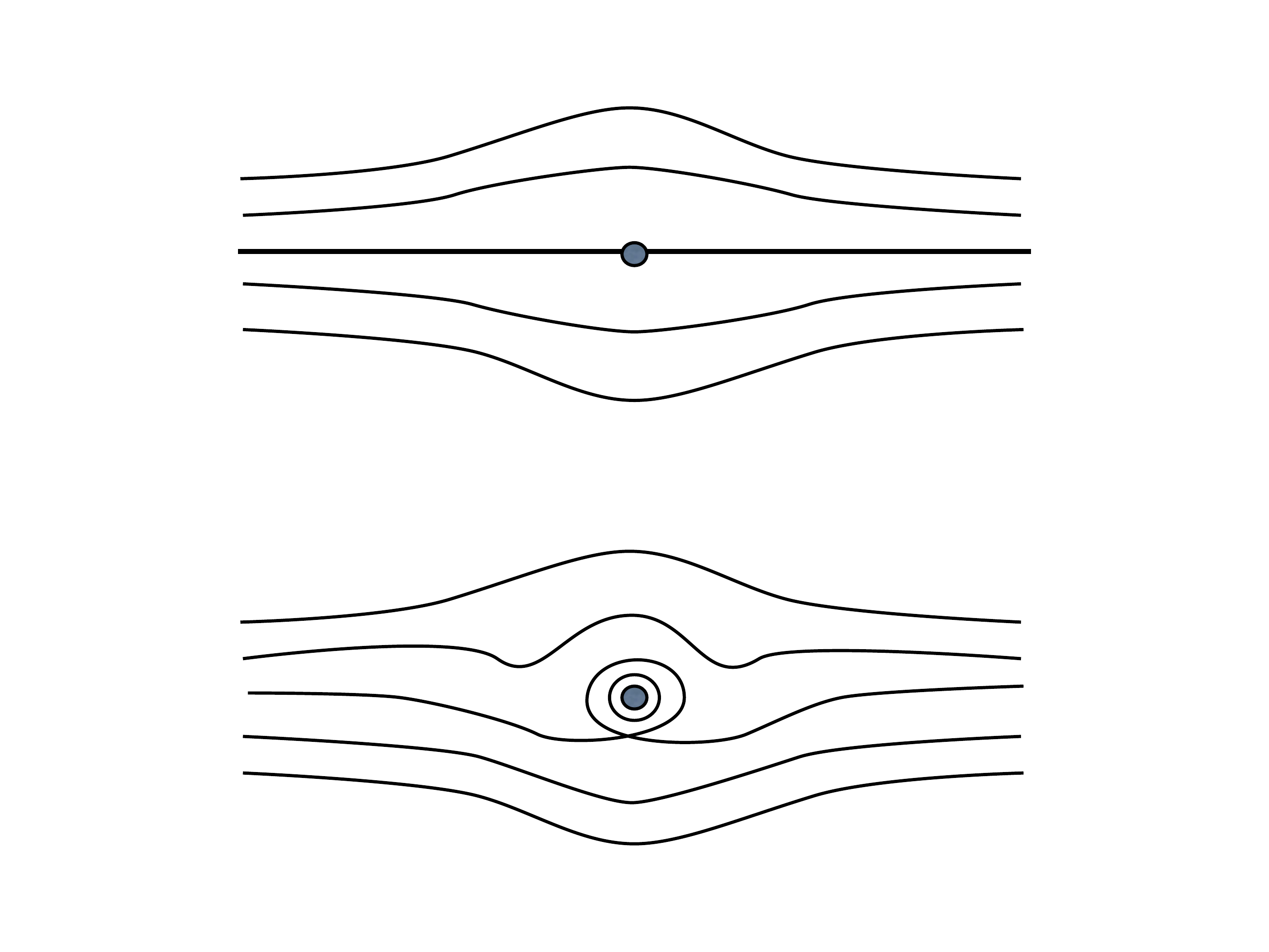}
\includegraphics[scale=.25]{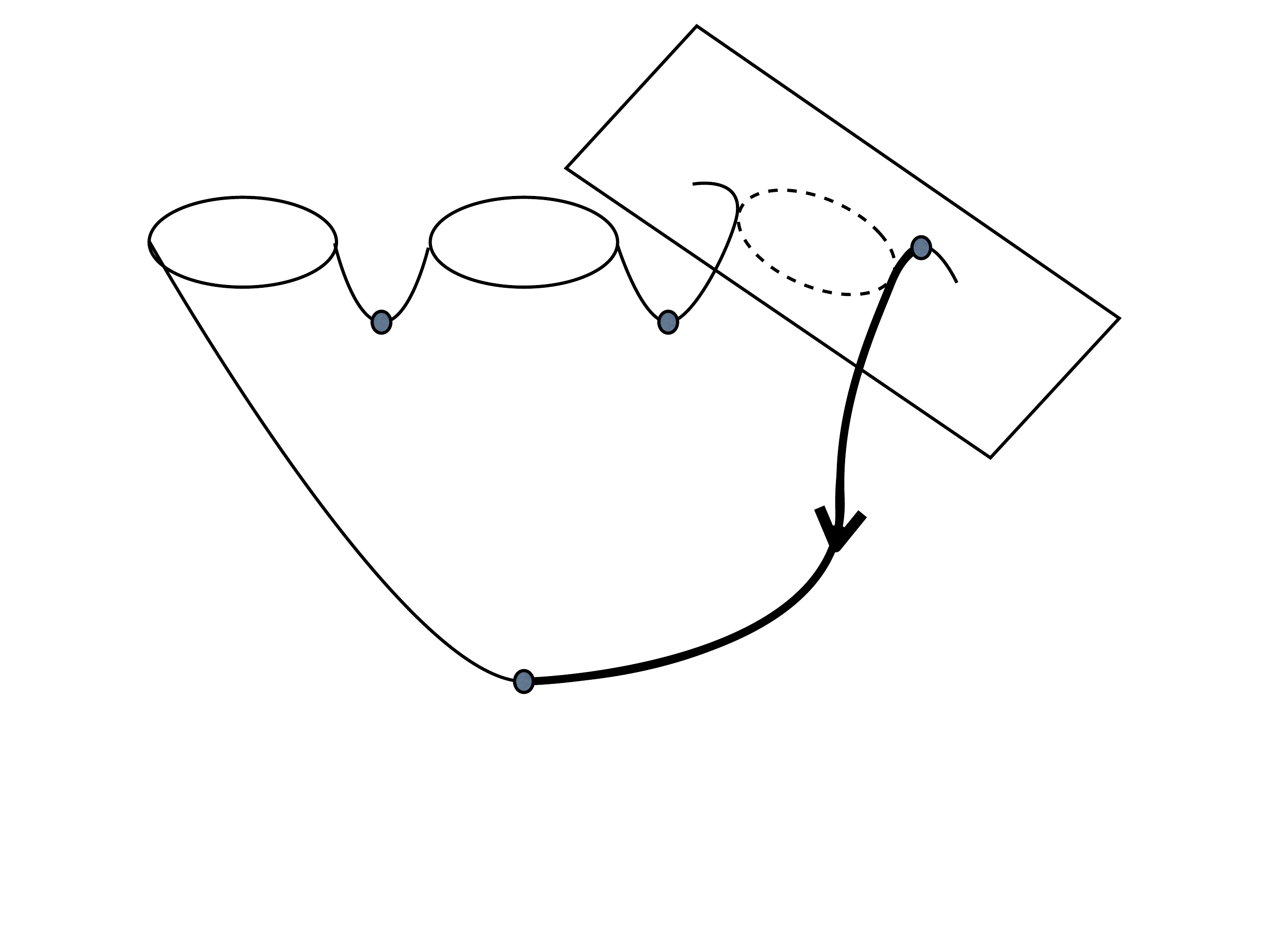}
\caption{\emph{Top}: Level sets of an unperturbed IIId boundary as a puncture. \emph{Middle}: Level sets of the perturbation. Note that the central rotating disk is excised to glue with other components. \emph{Bottom}: The unique flow-line.}
\label{fig1}
\end{center}
\end{figure}

We consider the one-boundary-component case first. For the purposes of computing rank, we may assume the extremum is a minimum by duality. Suppose first that we have a type IIId boundary component abutting our fixed component $S$. In \cite{c} we show that there is precisely one flow-line from each of the $p$ type IIId fixed points to the minimum, see Figure \ref{fig1}. In this case we have a cancellation because we are working with field coefficients and whatever element (even zero) of $N_h'/tors$ this flow-line corresponds to under the representation $\rho_{wm}$, it corresponds to a nonzero and thus invertible element of $Q(\ZZ/2[N_h'/tors])$. Thus the rank of the summand corresponding to the Nielsen class $\eta$ is $|\ind(\eta)|$.

Suppose next that $S$ does not abut any type IIId boundary components. Denote the minimum by $y$. Then for every saddle point $x$ we have two flow-lines to $y$.

\begin{figure}
\begin{center}
\includegraphics[scale=.25]{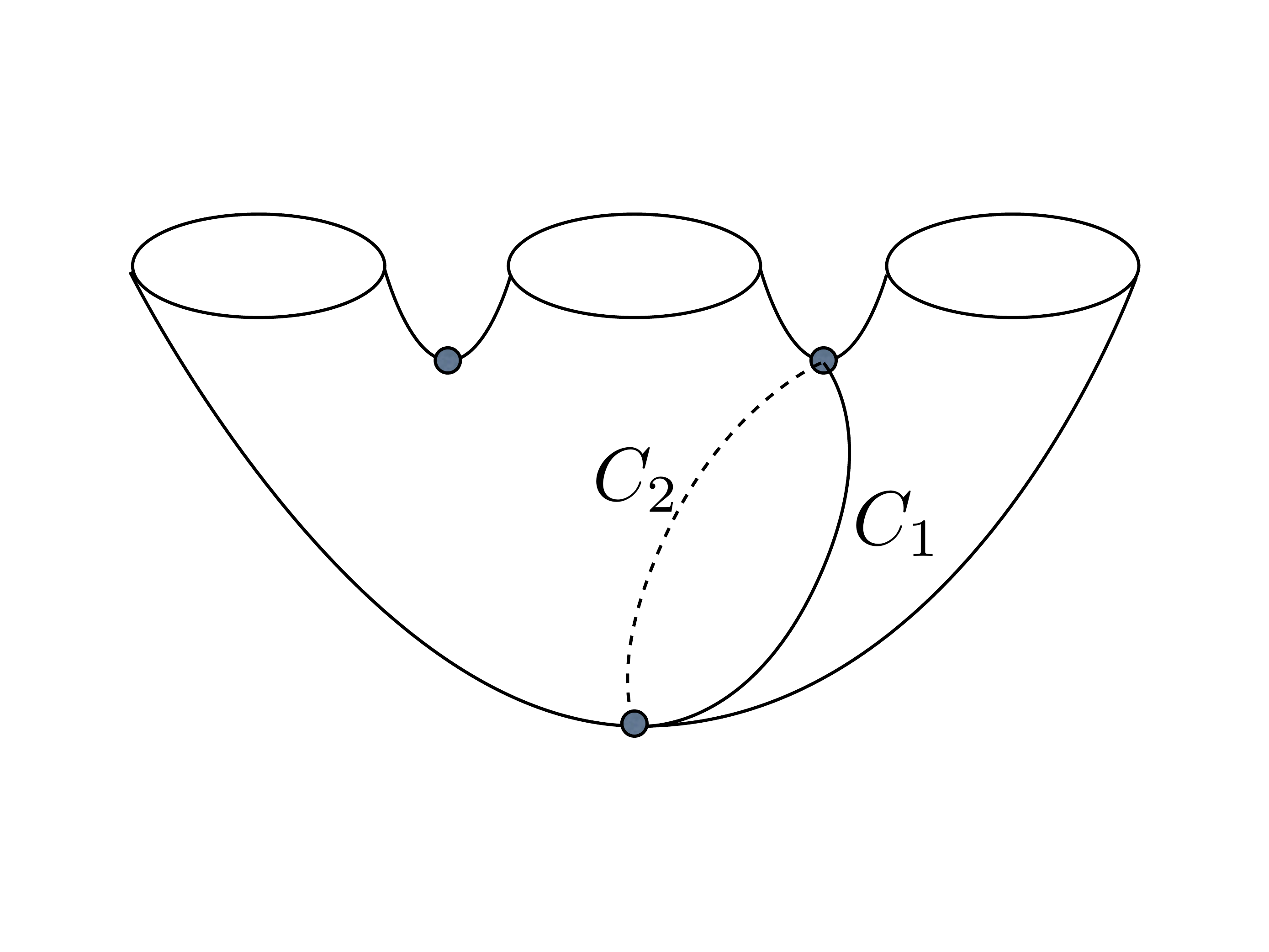}
\caption{The closure of the descending manifold.}
\label{fig2}
\end{center}
\end{figure}

\begin{lemma}
\label{desc}
We have that $\del x = a_x y$, where $a_x$ is nonzero if and only if the class of the closure of the descending manifold of $x$ in $N_h'/tors \subset H_1(\Sigma)$ is nonzero.
\end{lemma}

\proof
See Figure \ref{fig2}. We have two flow-lines $C_1$ and $C_2$. Thus $\del x = \rho([C_1\cdot C^{-1}_{y,x}]) y + \rho([C_2\cdot C^{-1}_{y,x}]) y = (a-b)y$. Here $a,b\in N_h'/tors$ are such that $a-b$ in $N_h'/tors$ is the class of the descending manifold of $x$. Note that the coefficient in $(a-b)y$ is not this same $a-b \in N_h'/tors$ but rather lies in the group ring (in which addition in $N_h'/tors$ would be written multiplicatively). In both situations, however, the result is zero if and only if $a = b$. The result follows.
\qed

If $a_x \neq 0$ for some $x$ then $y$ is a boundary because our coefficients lie in a field, and so we have a cancellation and the rank is only $|\ind(\eta)|$. If $a_x = 0$ for all $x$, this is the statement that the differential vanishes in the summand corresponding to the Nielsen class $\eta$, and so the rank is $|\ind(\eta)|+2$.

\begin{figure}
\begin{center}
\includegraphics[scale=.25]{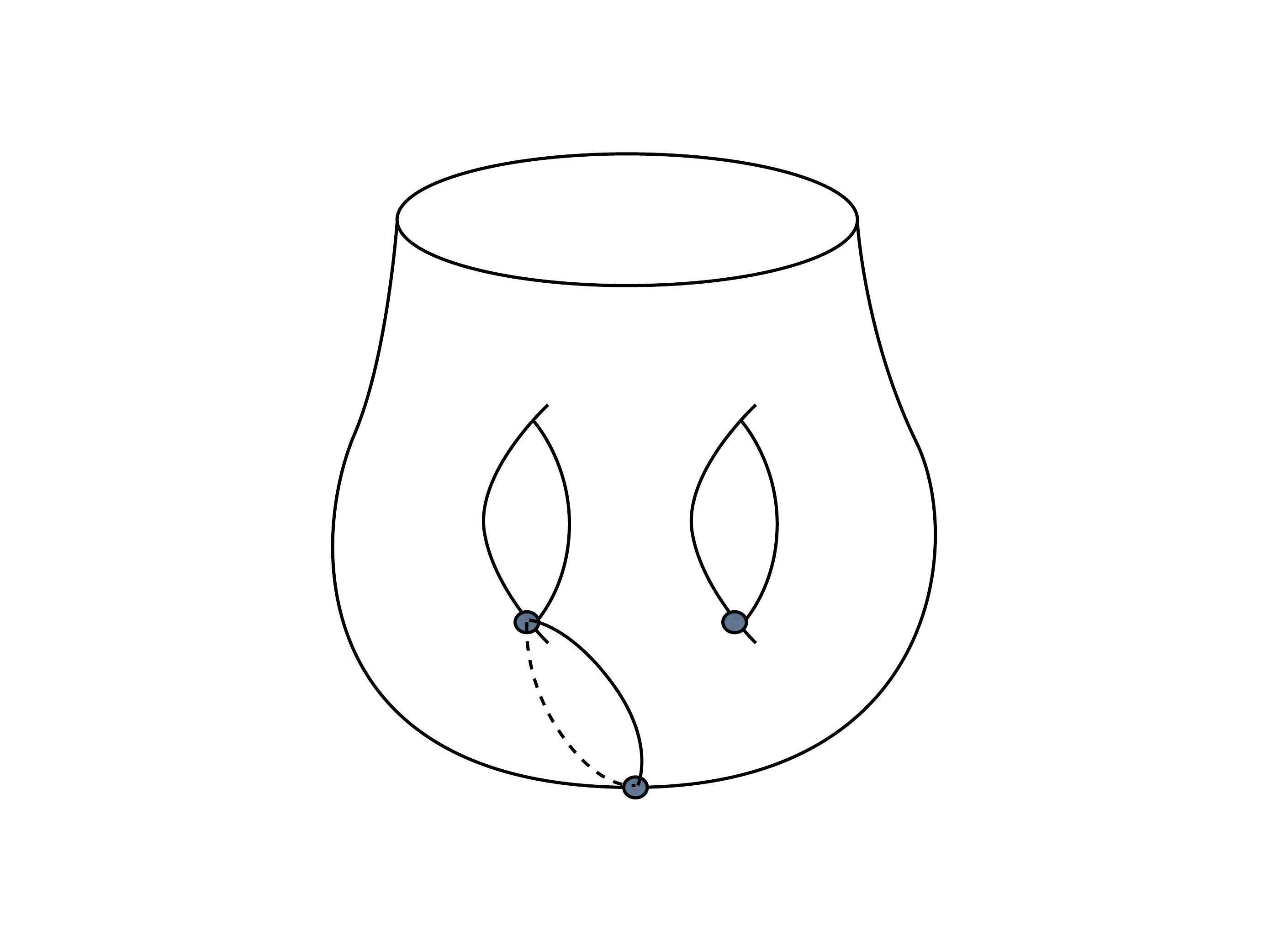}
\caption{The nonzero genus case.}
\label{fig3}
\end{center}
\end{figure}

\begin{figure}
\begin{center}
\includegraphics[scale=.35]{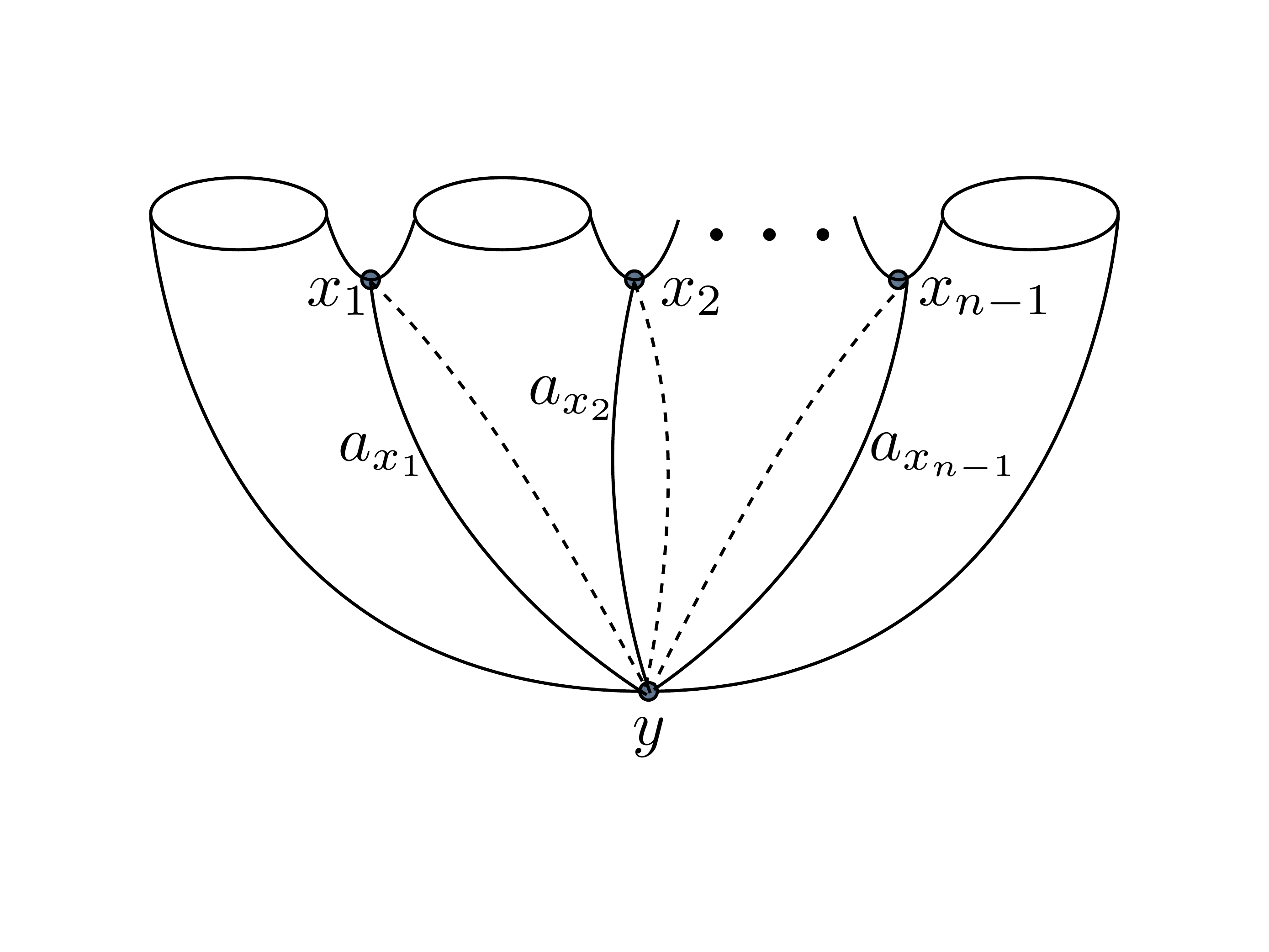}
\caption{The genus zero case with all boundary components rotating the same direction.}
\label{fig4}
\end{center}
\end{figure}

\begin{lemma}
In the above situation, $a_x = 0$ for every saddle point $x$ if and only if $S$ has genus zero and every boundary component of $S$ is nulhomologous in $H_1(\Sigma,\del\Sigma)$.
\end{lemma}

\proof
In the case in which $S$ has genus, we can find two descending manifolds which meet algebraically once. It follows that neither of them are nulhomologous. See Figure \ref{fig3}.

If $S$ is genus zero, it appears as in Figure \ref{fig4}. The homology class of each of the boundary components is given by either (plus or minus) the homology class of one of the descending manifolds (if on either end) or the difference of two such. We see all the descending manifolds are nulhomologous if and only if all of the boundary components are.
\qed

In the case in which $S$ has no boundary, we see that $a_x$ is never zero for saddle points $x$ and thus we cancel the minimum. Dually, we also cancel the maximum, say with any saddle point we haven't used to cancel the minimum. Thus in this case the rank is $|\ind(\eta)|$. Summing up, we've shown:

\begin{prop}
\label{2extra}
Consider a Nielsen class $\eta$ corresponding to a fixed component $S$. If $S$ does not abut any type IIId boundary components, every boundary component rotates in the same direction, $S$ has genus zero, and every boundary component of $S$ is nulhomologous in $H_1(\Sigma,\del\Sigma)$ then the rank of the summand of $HF_*(\phi,Q(\ZZ/2[N_h'/tors]))$ corresponding to such a Nielsen class $\eta$ is $|\ind(\eta)|+2$. Otherwise it is $|\ind(\eta)|$.\qed
\end{prop}

\subsection{Construction of maps using symplectic vector fields}

We now construct a map $\psi$ in our mapping class $h$ such that in every Nielsen class $\eta$, the number of fixed points in this Nielsen class equals the rank of the summand of $HF_*(\phi,Q(\ZZ/2[N_h'/tors]))$ corresponding to the Nielsen class $\eta$, which we computed in the previous section.

We start with a map $\phi_{st}$, a standard-form reducible map but which is the identity on any fixed components as opposed to having been perturbed by Hamiltonian flows. We next perturb by Hamiltonian flows on any components for which there was no cancellation in the previous section; that is, components for which the number of fixed points equalled the rank. These are components which have boundary components rotating in different directions as well as components which satisfy all of the criteria in Propostion \ref{2extra}. 

Next we use a modified Hamiltonian on components $S$ which meet a type IIId boundary, geometrically cancelling the extremum with one of the type IIId fixed points. We start with a Hamiltonian function with at most one extremum on such a component which patches together with the Hamiltonian on a neighborhood of the boundary of the pseudo-Anosov region as in the previous section.

\begin{figure}
\begin{center}
\includegraphics[scale=.25]{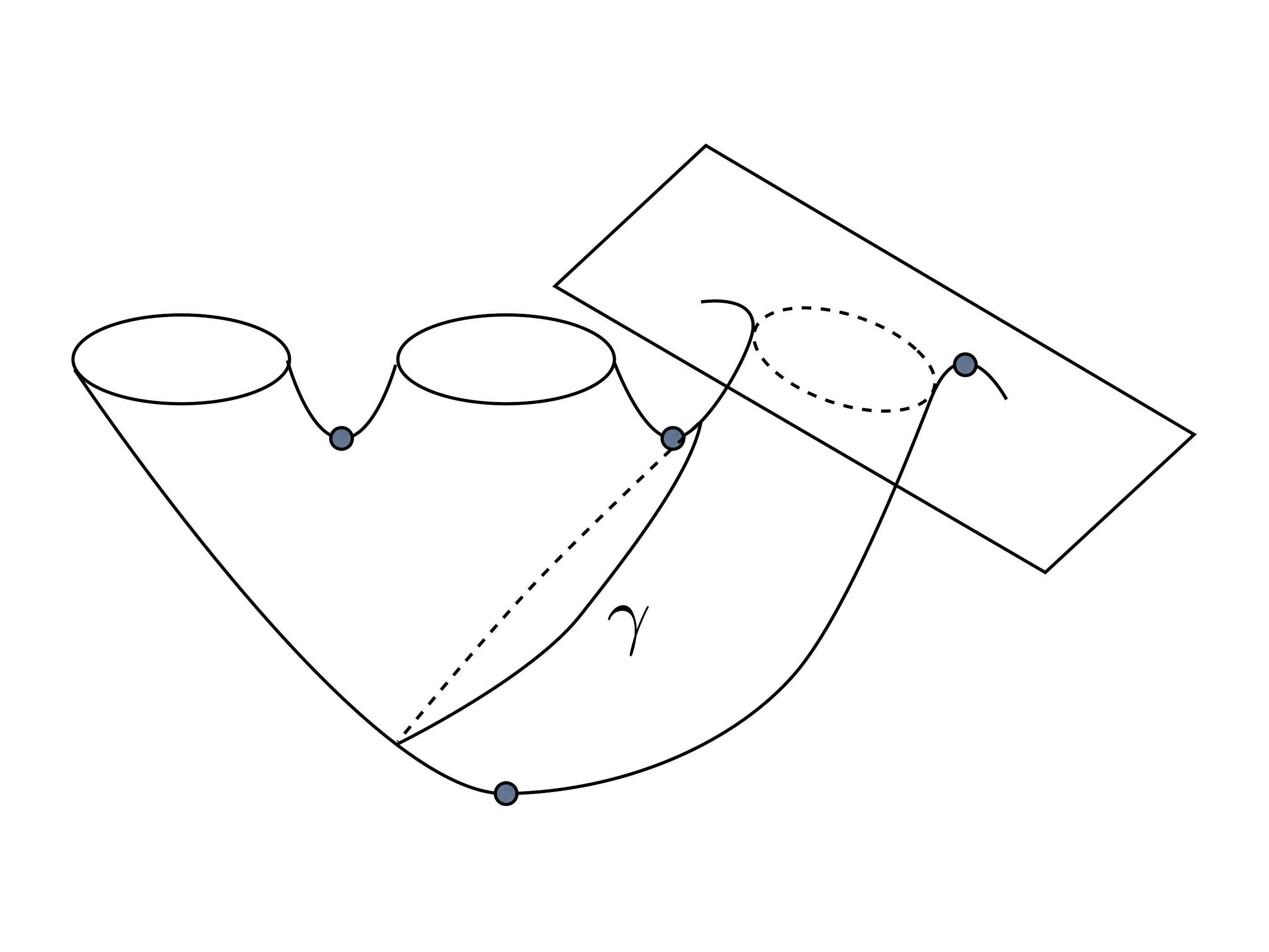}Ê
\includegraphics[scale=.25]{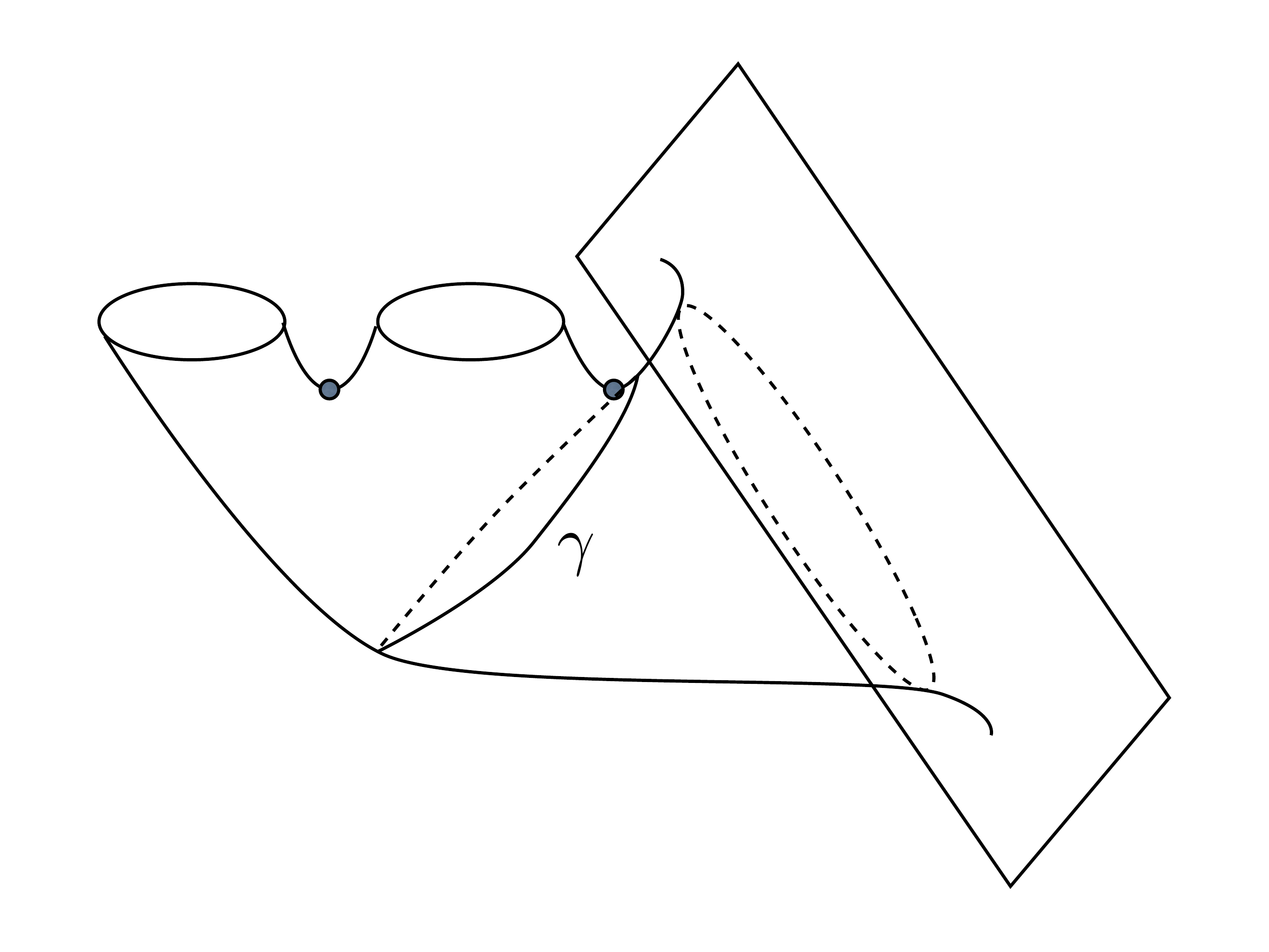}Ê
\caption{Geometrically cancelling the minimum when meeting a type IIId boundary. \emph{Top}: Before. \emph{Bottom}: After.}
\label{fig5}
\end{center}
\end{figure}

\begin{figure}
\begin{center}
\includegraphics[scale=.25]{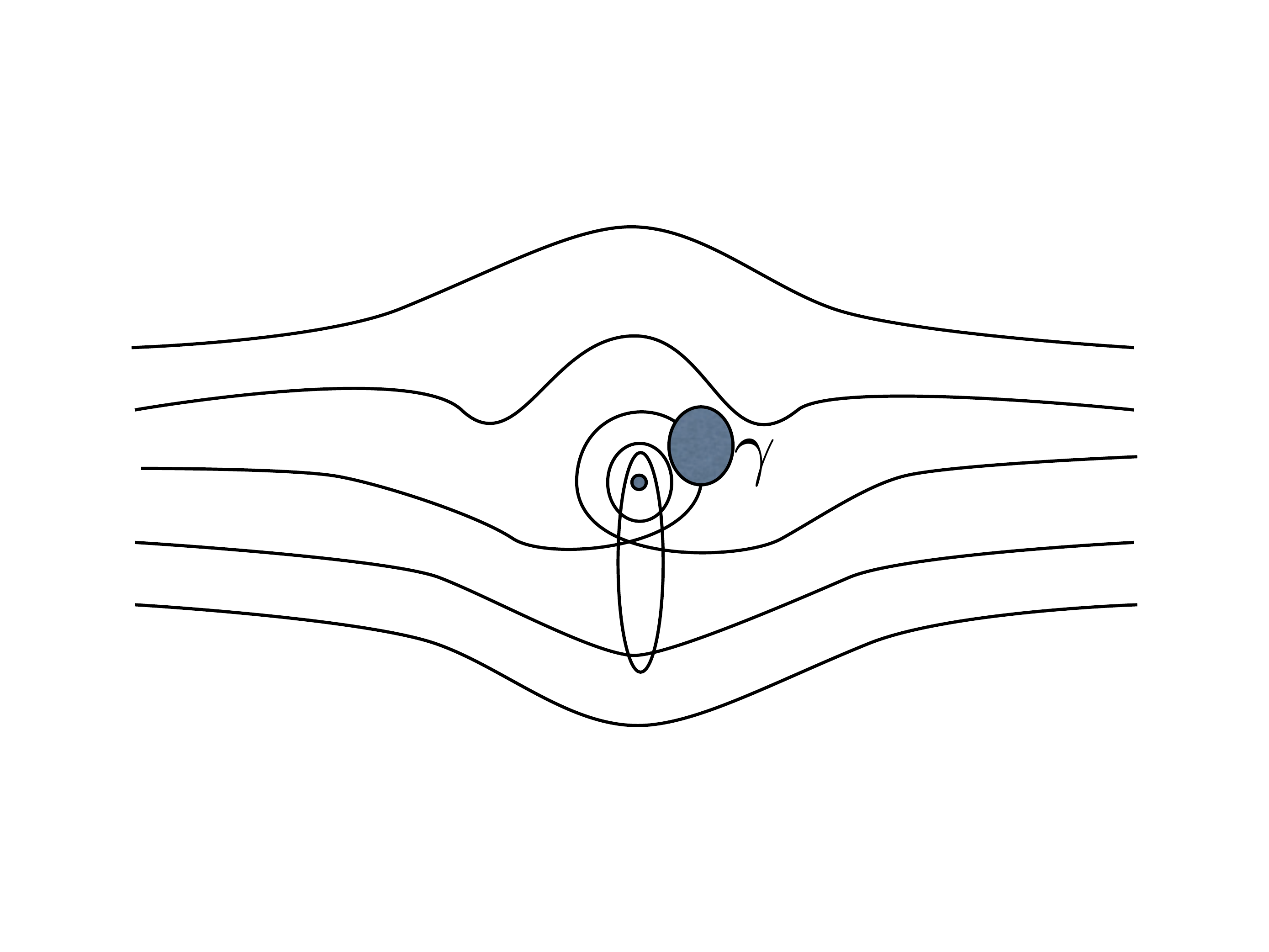}
\includegraphics[scale=.25]{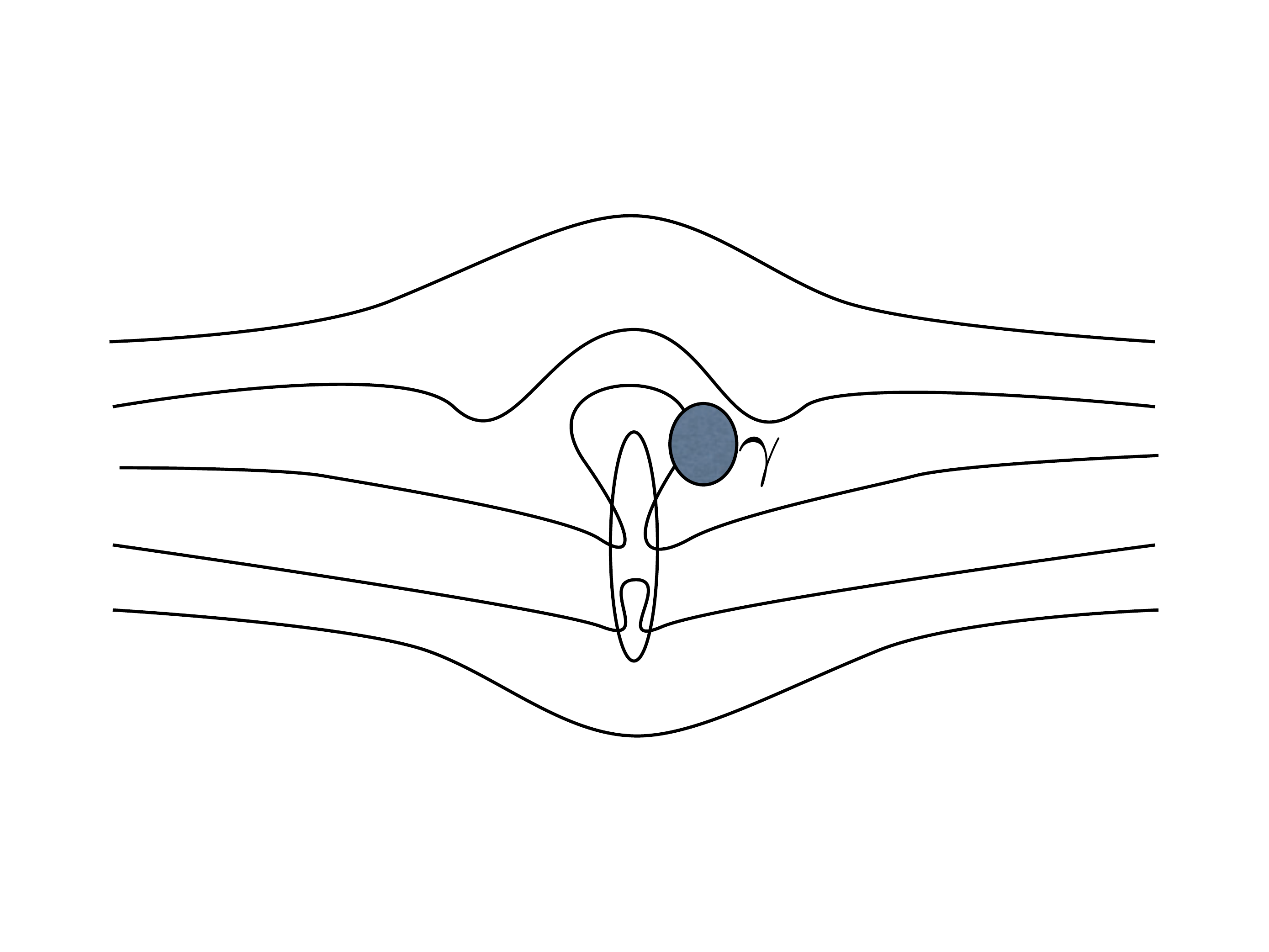}
\caption{The perturbed type IIId boundary with shaded $\gamma$-bounded disk to be excised and replaced with the rest of the component $S$, and unshaded region in which we are modifying the Hamiltonian. \emph{Top}: Before. \emph{Bottom}: After.}
\label{fig6}
\end{center}
\end{figure}

\begin{lemma}
On a component $S$ meeting a type IIId boundary, there exists a modification of the with the aforementioned Hamiltonian perturbation whose critical points are all nondegenerate saddle points.
\end{lemma}

\proof
We geometrically cancel one of the $p$ type IIId fixed points and the fixed point corresponding to the minimum as in Figure \ref{fig5}. To do this, we consider, as in Figure \ref{fig6}, the perturbed situation but without the central disk excised. We draw a loop $\gamma$, excising the disk it bounds and replacing it with ``the rest'' of $S$, i.e. the portion to the left of $\gamma$ in Figure \ref{fig5}. This reduces the situation to one standard case. We rescale the Hamiltonian on $S$ to be small if necessary so there are points on which the Hamiltonian evaluates to a number less than the evaluation of the Hamiltonian at the (local) minimum. Then there is no problem.
\qed

Next we consider components $S$ whose boundary components rotate all in the same direction and which meet no type IIId boundaries, but have nonzero genus.

\begin{figure}
\begin{center}
\includegraphics[scale=.25]{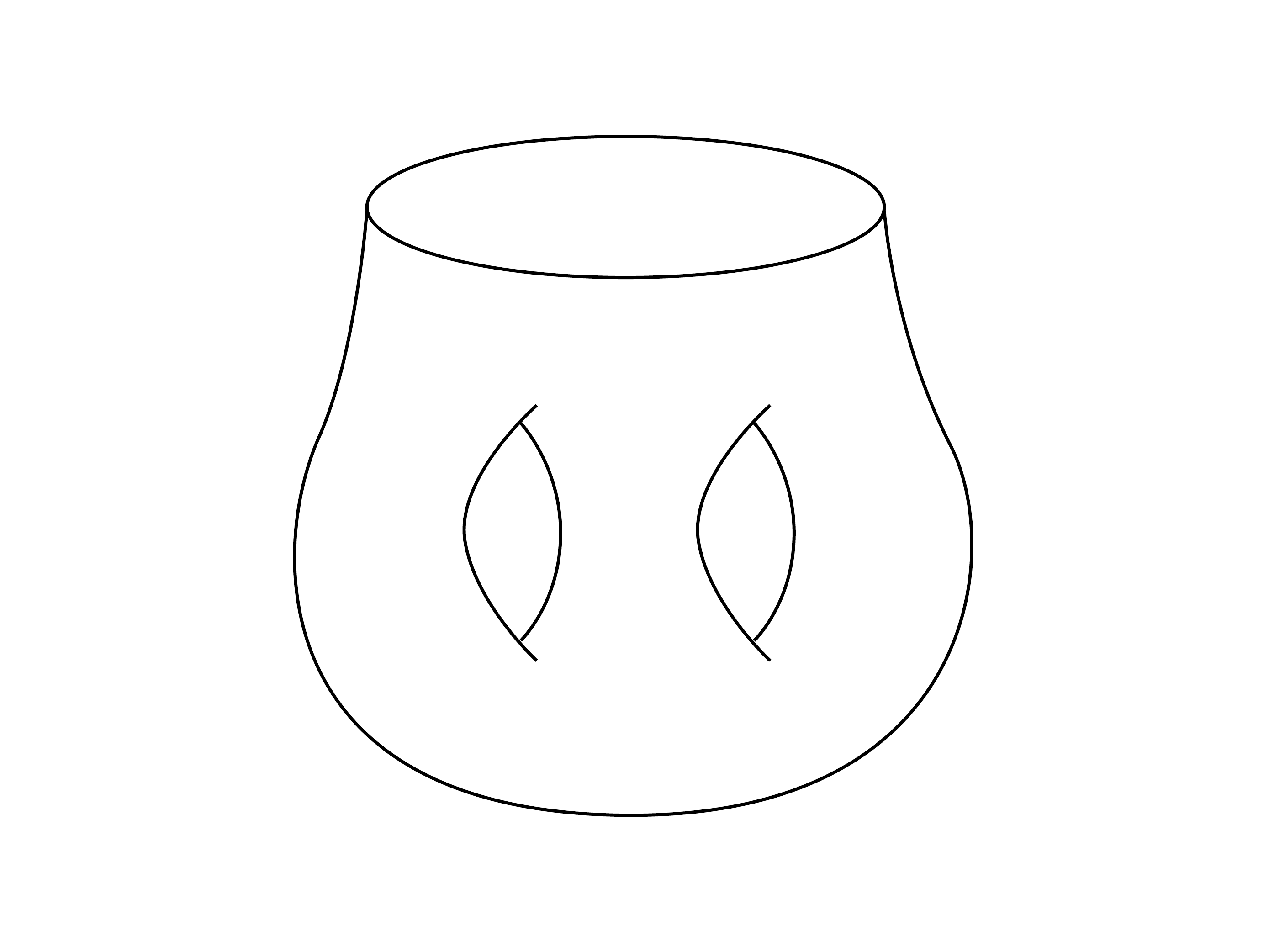}
\includegraphics[scale=.25]{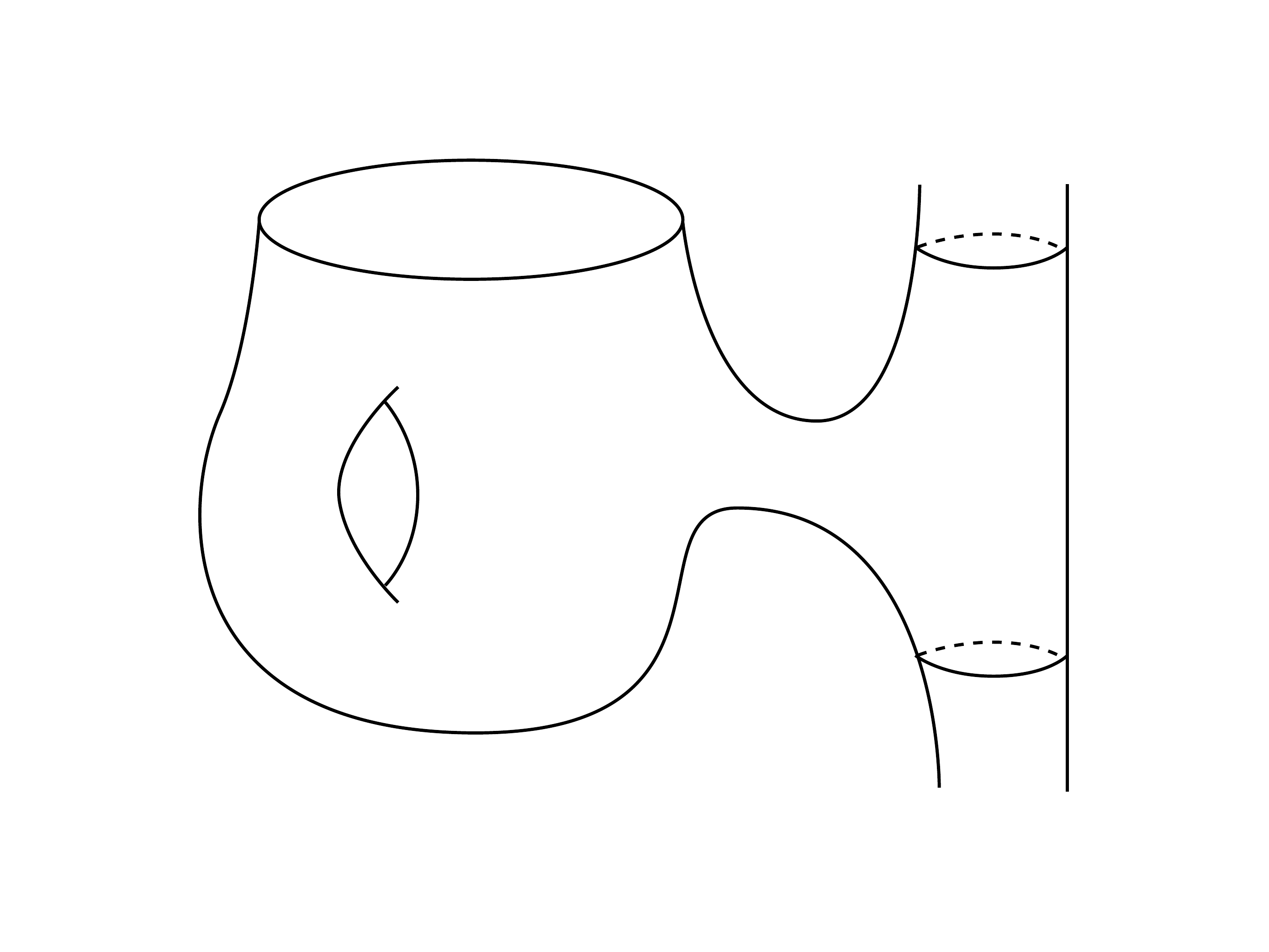}
\caption{On a component which has genus, we can locally cancel the minimum geometrically with an $S^1$-valued Hamiltonian. \emph{Top}: Before. \emph{Bottom}: After.}
\label{fig7}
\end{center}
\end{figure}

\begin{lemma}
There exists an $S^1$-valued Hamiltonian on such a component whose associated symplectic vector field is parallel to the boundary and whose critical points are all nondegenerate saddle points.
\end{lemma}

\proof
We modify the Hamiltonian away from the boundary as in Figure \ref{fig7}.
\qed

Now we perturb this whole map by the flow of a small symplectic vector field.

\begin{lemma}
There exists a symplectic vector field $V$ transverse to every reducing curve $C$ for which $[C] \neq 0$ in $H_1(\Sigma,\del\Sigma)$.
\end{lemma}

\proof
The plan is to choose a totally irrational Flux class. In order to see that we can choose the flow to be transverse to every reducing curve, we choose a handlebody bounding $\Sigma$ such that the reducing curves all bound disks. Next we collapse this to the underlying graph and put flows on each edge of the graph whose only rational dependances are given by balancing conditions at the vertices. We then use this as a guide to build the flow on the surface.
\qed

We rescale the symplectic vector field so that it is small enough that its time-$1$ flow does not move any of the fixed points too much and does not create any new fixed points.

For components $S$ whose boundary components rotate all in the same direction, which meet no type IIId boundaries, and are genus zero but which have a non-nulhomologous boundary component, we now modify the symplectic vector field $V$ in a neighborhood of the component, keeping it transverse to each reducing curve.

\begin{figure}
\begin{center}
\includegraphics[scale=.25]{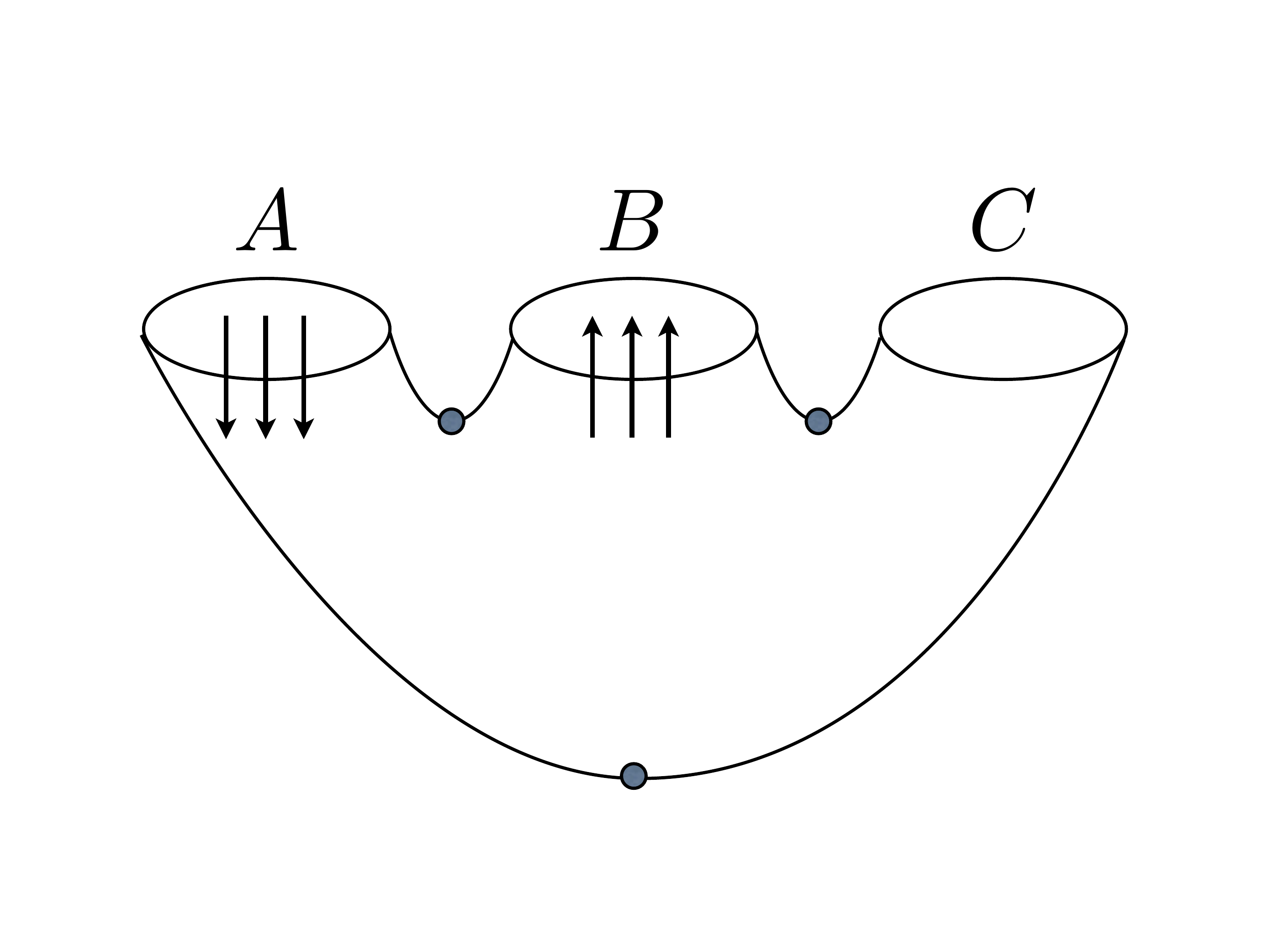}
\includegraphics[scale=.25]{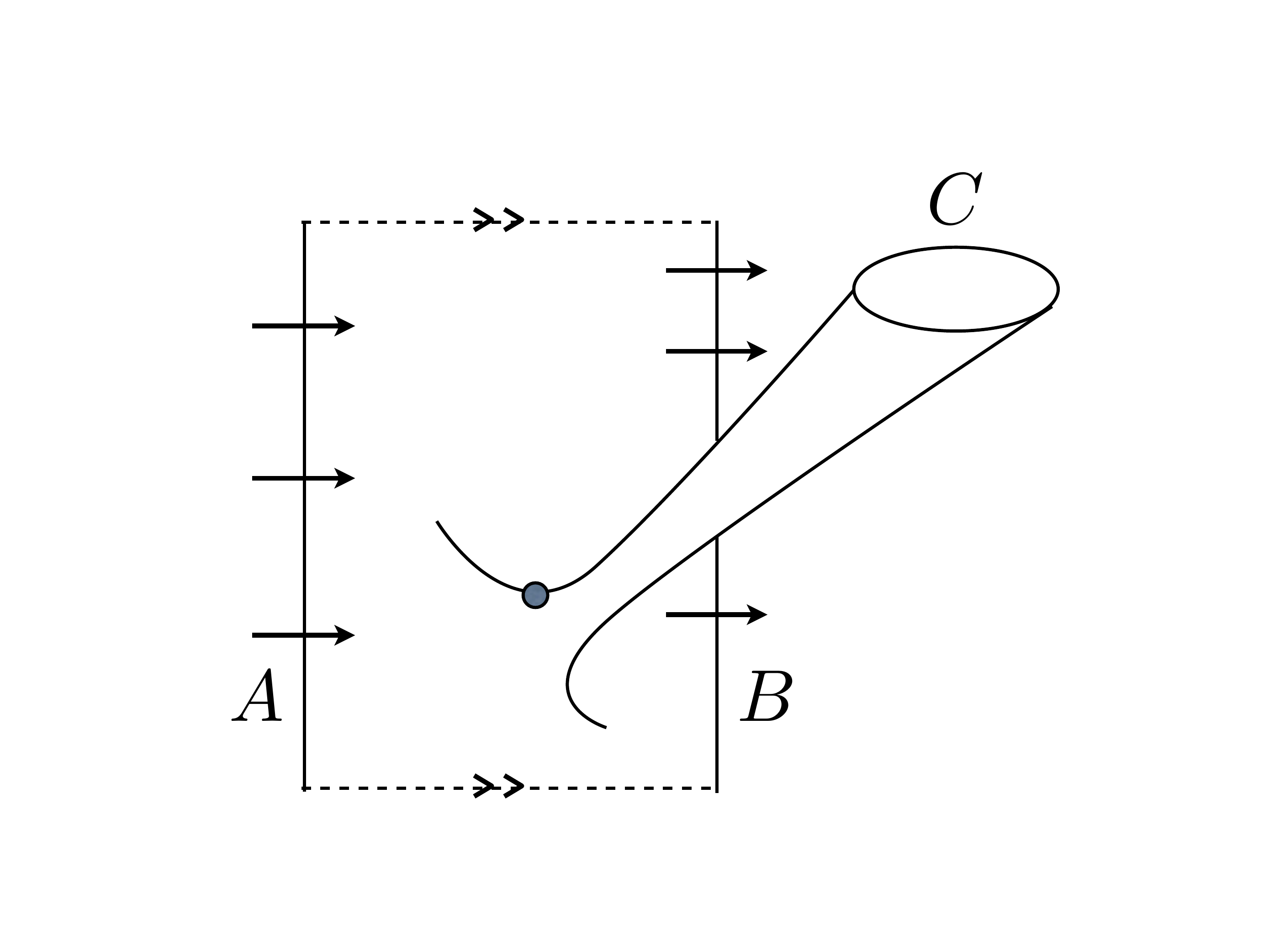}
\caption{On a genus zero component, if the symplectic vector field is nonzero on at least one boundary component, then we can geometrically cancel the minimum with an $S^1$-valued Hamiltonian.  \emph{Top}: Before. \emph{Bottom}: After.}
\label{fig9}
\end{center}
\end{figure}

\begin{lemma}
There exists an $S^1$-valued Hamiltonian on such a component whose associated symplectic vector field is transverse to the boundary and whose critical points are all nondegenerate saddle points.
\end{lemma}

\proof
We modify the Hamiltonian away from the boundary as in Figure \ref{fig9}.
\qed

We've shown:

\begin{prop}
Given a reducible mapping class $h$, there exists a map $\psi$ in the mapping class such that in every Nielsen class $\eta$, the number of fixed points in this Nielsen class equals the rank of the summand of $HF_*(\phi,Q(\ZZ/2[N_h'/tors]))$ corresponding to the Nielsen class $\eta$.\qed
\end{prop}

Combining this with our observations regarding identity, periodic, pseudo-Anosov mapping classes and, for reducible mapping classes, with Proposition \ref{2extra}, we obtain (summing over Nielsen classes):

\begin{thm}
\label{fixedboundsection}
The minimum number of fixed points of an area-preserving map $\phi$ with nondegenerate fixed points in a mapping class $h$ is given by:
$$ \left\{ \begin{array}{cc}
\sum_\eta \left| \ind(\eta) \right| & h\ \textrm{is periodic or pseudo-Anosov} \\
\sum_\eta \left| \ind(\eta) \right| + 2 A & h\ \textrm{is reducible} \end{array} \right.$$
where $A$ is the number of genus zero components of the reducible mapping class on which the map is the identity, which do not abut any pseudo-Anosov components, and all of whose boundary components rotate in the same direction and are nulhomologous or homologically boundary parallel.\qed
\end{thm}

\section{Degenerate fixed points}
\label{degenerate}

If we are allowed degenerate fixed points, each of the Nielsen classes $\eta$ for which the bound for nondegenerate fixed points was $|\ind(\eta)|$ can be reduced to a single degenerate fixed point. To see this, note that the only of these situations in which we are not already reduced to a single fixed point are in cases in which the Nielsen class is associated with a $p$-prong pseudo-Anosov singularity or in which the Nielsen class is associated to a fixed component (possibly meeting type IIId boundaries). In the former case, we simply modify the singular Hamiltonian $H_{sing} = \mu r^2 \cos(p\theta) = \mu \Re(z^p)/|z|^{p-2}$ to a smooth Hamiltonian which agrees with $H_{sing}$ outside a small ball and inside a yet smaller ball is $C r^p \cos(p\theta) = C \Re(z^p)$, which is smooth at the origin, where it has a (generalized) monkey saddle. In the latter case, all of the fixed points are index negative one, i.e. are given by saddles, and we can again combine them all into the appropriate generalized monkey saddle, using $p = |\ind(\eta)|+1$ (the index of such a degenerate fixed point is $1-p$).

\begin{figure}
\begin{center}
\includegraphics[scale=.25]{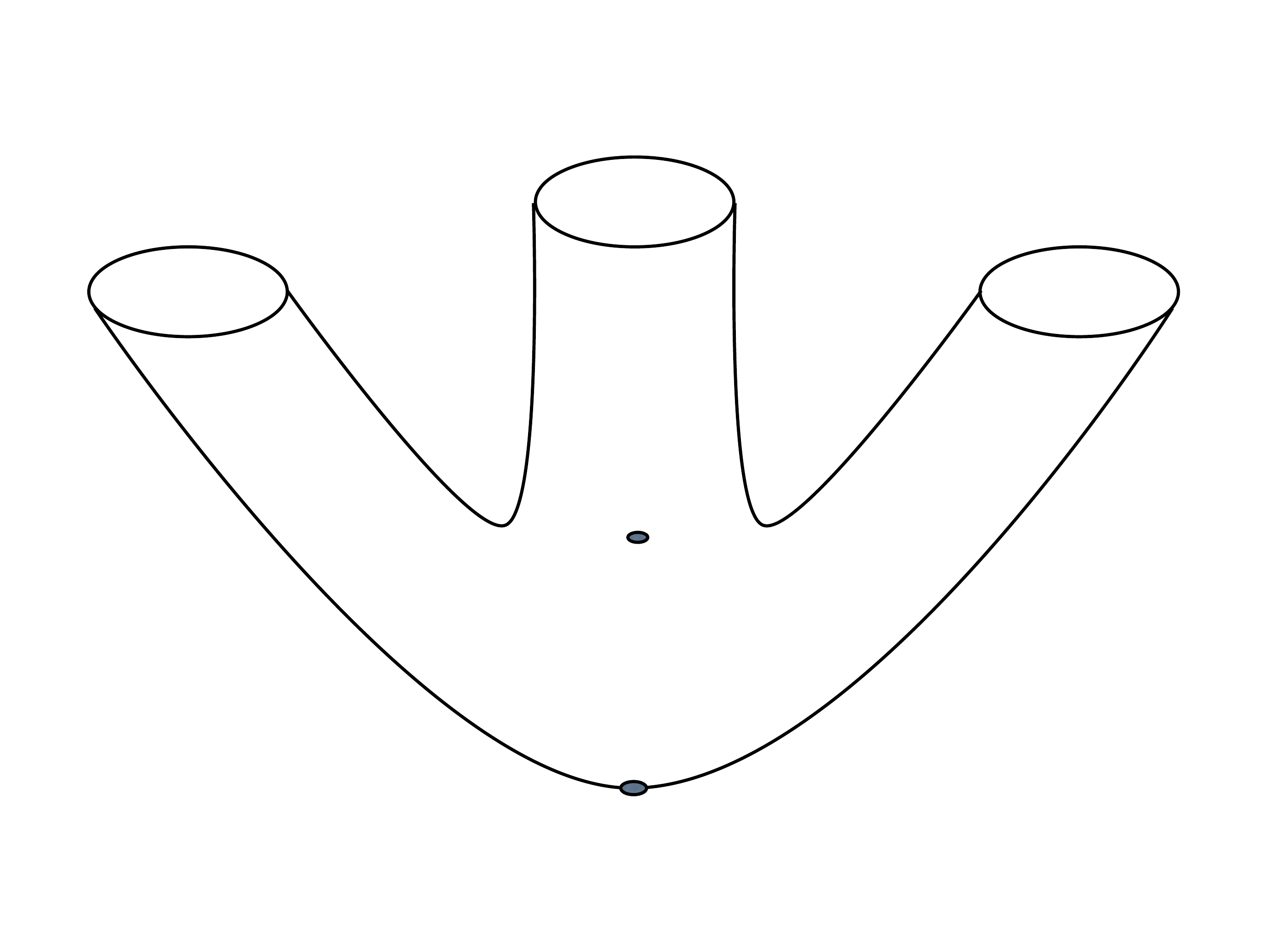}
\caption{A three boundary component genus zero component with one minimum and, instead of two nondegenerate index negative one fixed points, one $3$-prong monkey saddle.}
\label{fig10}
\end{center}
\end{figure}

Similarly, each of the Nielsen classes $\eta$ for which the bound for nondegenerate fixed points was $|\ind(\eta)|+2$ can be reduced to two fixed points: each of these corresponds to a genus zero component of the standard form map on which the map is the identity and for which each boundary component rotates in the same direction. Each of these in the nondegenerate case one index one fixed point, some number $k$ of index negative one fixed points, where the component has $k+1$ boundary components. All of the index negative one fixed points can be combined in one degenerate fixed point of the same sort as in the previous paragraph. See Figure \ref{fig10}. Our task is now to show that we can do no better. That is, we cannot combine these two fixed points into one.

To show this, we use a cohomology operation. The argument has similarities with arguments that cup lengths give bounds on (even degenerate) fixed points, but even though we have, by \cite[\S 2.5]{c}, a deformation-invariant module structure over the quantum homology of $\Sigma$ (which agrees with $H_*(\Sigma)$), the module structure is trivial. Every element of $H_1(\Sigma)$ acts as zero because each of the descending manifolds from Lemma \ref{desc} is nulhomologous, and thus has zero algebraic intersection with any element of $H_1(\Sigma)$.

All is not lost, however. If we restrict our attention to one Nielsen class $\eta$, there is a sense in which this descending manifold is homologically essential for this Nielsen class.

\begin{lemma}
\label{pi1isS}
Let $\Sigma$ be a surface of negative Euler characteristic, $h$ reducible mapping class, and $\eta$ a Nielsen class corresponding to a fixed component $S$. Then $H_1(\Gamma(M_\phi)_\eta) \cong H_1(S)$ for any $\phi \in h$. In fact, $\pi_1(\Gamma(M_\phi)_\eta) \cong \pi_1(S)$. Furthermore, the image of $$H_1(\Gamma(M_\phi)_\eta) \rightarrow H_2(M_\phi) \rightarrow H_1(\Sigma)$$ in $H_1(\Sigma)$ agrees with the image of $H_1(S) \rightarrow H_1(\Sigma)$.
\end{lemma}

\proof
We consider $\gamma$ as a map $\RR\rightarrow\Sigma$ with $\gamma(t) = \phi(\gamma(t+1))$. An element of $\pi_1(\Gamma(M_\phi),\gamma)$ is of the form $\gamma_s(t)$ for $s\in S^1 = \RR/\ZZ$ with $\gamma_0(t) = \gamma(t)$.
We consider $\alpha_0(s) = \gamma_s(0)$ and $\alpha_1(s) = \gamma_s(1)$. These are closed curves on $\Sigma$ and $\phi(\alpha_1(s)) = \alpha_0(s)$. Furthermore, $\alpha_0(s)$ is homotopic to $\alpha_1(s)$ by the homotopy $\alpha_t(s) = \gamma_s(t)$.

As in \cite[Lemma 3.2]{c}, we have a fibration $\Omega\Sigma \rightarrow \Gamma(M_\phi) \rightarrow \Sigma$ and thus a long exact sequence on homotopy groups, a piece of which is \beqa \pi_2(\Sigma,\gamma(0)) \rightarrow \pi_1(\Gamma(M_\phi), \gamma) \rightarrow \pi_1(\Sigma,\gamma(0)).\eeqa The image in $\pi_1(\Sigma,\gamma(0))$ of the element of $\pi_1(\Gamma(M_\phi),\gamma)$ represented by the homotopy $\gamma_s(t)$ is represented by $\alpha_0(s)$. Because $\Sigma$ is a surface of negative Euler characteristic, $\pi_2(\Sigma) = 0$. Thus we have an injection \beqa \pi_1(\Gamma(M_\phi), \gamma) \hookrightarrow \pi_1(\Sigma,\gamma(0)).\eeqa

We choose our map $\phi$ (amongst those in the mapping class) to be a standard form one which is the identity on $S$ and choose our basepoint $\gamma$ to be the constant path at some point in $S$. We claim that the image of this map is $\pi_1(S,\gamma(0))$. We note that the image consists of elements of $\pi_1(\Sigma,\gamma(0))$ represented by loops $\alpha(s)$ based at $\gamma(0)$ which are homotopic through loops based at $\gamma(0)$ (because $\gamma$ is the constant path) to $\phi(\alpha(s))$.Thus the image contains $\pi_1(S,\gamma(0))$.

We now claim that any $\alpha(s)$ based at $\gamma(0)$ homotopic to $\phi(\alpha(s))$ through loops based at $\gamma(0)$ can be homotoped inside $S$. This would give the result. This follows from \cite[Lemma 3.4]{jg}, with its modification to the standard form maps for the area-preserving case (in which we need to consider multiple parallel Dehn twist regions with fixed annuli in between) given in \cite[Lemma 3.8 and Corollary 3.9]{c}. These state that any path between two fixed points of a standard form map which is homotopic rel endpoints to $\phi$ applied to itself can be homotoped inside the fixed point set of $\phi$. The component of $\mathrm{Fix}(\phi)$ containing $\gamma(0)$ is simply $S$.

Finally, we note that the result continues to hold for any other map in the mapping class. Nielsen classes are well-defined on the entire mapping class $h$ by \cite[Lemma 4.2]{c}, so this statement is sensible.
\qed

Now we restrict our attention to fixed genus zero components $S$ which do not meet type IIId boundaries, and whose boundary components all rotate in the same direction and are all nulhomologous.

\begin{lemma}
\label{allokS}
The Floer homology chain complex $CF_*(\phi,\ZZ/2[H_1(\Gamma(M_\phi)_\eta)];\eta)$ of a map $\phi \in h$, restricted to Nielsen class $\eta$ summand of such a fixed component $S$, with coefficients in the representation $\ZZ/2[H_1(\Gamma(M_\phi)_\eta)]$, is well defined on the entire mapping class $h$ and invariant up to quasi-isomorphism. The same holds for $CF_*(\phi,\ZZ/2;\eta)$, where the coefficients are the trivial representation $\ZZ/2$.
\end{lemma}

\proof
It is well-defined for any $\phi$ because $\omega_\phi: H_1(\Gamma(M_\phi)_\eta) \rightarrow \RR$ is the zero map for any $\phi$ (which implies that every $\phi$ is $\eta$-weakly monotone as defined in \cite{c}). This follows because $\omega_\phi: H_1(\Gamma(M_\phi)_\eta) \rightarrow \RR$ agrees with the map $$H_1(\Gamma(M_\phi)_\eta) \longrightarrow \left(\ker(c_\phi) \subset H_2(M_\phi)\right) \longrightarrow \left(N_h \subset H_1(\Sigma)\right) \stackrel{\Flux}{\longrightarrow} \RR$$ and the image of $H_1(\Gamma(M_\phi)_\eta)$ under all but the last composition is the image of $H_1(S) \rightarrow H_1(\Sigma)$, by Lemma \ref{pi1isS}. This, however, is zero because $S$ has genus zero so that $H_1(\del S)$ surjects to $H_1(S)$, but we've assumed the boundary of $S$ is nulhomologous in $\Sigma$. The invariance up to quasi-isomorphism follows from Theorem \ref{invariance}, part (2).
\qed

We compute with $\phi$ a standard form map. We assume without loss of generality that the Morse-Smale function on the component $S$ has one extremum, a minimum (otherwise reverse orientation on $\Sigma$). The homology $HF_*(\phi,\ZZ/2;\eta)$ has dimension $|\ind(\eta)|+2$, generated by a fixed point $y$ corresponding to the minimum and $|\ind(\eta)|+1$ fixed points $x_i$ of index negative one.

Consider a homomorphism \beqa \beta \in &&\Hom(H_1(\Gamma(M_\phi)_\eta),\ZZ/2) = \Hom(H_1(S),\ZZ/2) \\ &&= H^1(S;\ZZ/2) = H_1(S,\del S;\ZZ/2).\eeqa As in \cite[\S 12.1.3]{hs2} and \cite[\S 2.5]{c} we get a degree one map $$\del_\beta : HF_*(\phi,\ZZ/2;\eta) \rightarrow HF_*(\phi,\ZZ/2;\eta)$$ defined by $$\del_\beta z = \sum_w \sum_{C \in \MM_1(z,w)} \beta([C\cdot C^{-1}_{z,w}])\cdot w.$$ Moreover, we show in \cite[Proof of Proposition 2.9]{c} that $\del_\beta$ is well-defined (purely algebraically) up to quasi-isomorphism of the pair $$CF_*(\phi,\ZZ/2[H_1(\Gamma(M_\phi)_\eta)];\eta) \ \textrm{and} \ CF_*(\phi,\ZZ/2;\eta).$$ Thus by Lemma \ref{allokS} $\del_\beta$ is a well-defined operation $HF_*(h,\ZZ/2;\eta)$ on the $\eta$ component of the Floer homology for the mapping class $h$. This operation is what replaces quantum cup products such as in \cite{sc2} in a cup-length argument.

\begin{lemma}
\label{beta}
There is a $\beta$ such that $\del_\beta x_i = y$. In fact, for any sum $\sum_i c_i x_i$ with $c_1 \in \ZZ/2$ not all zero, there exists a $\beta$ such that $\del_\beta \sum_i c_i x_i = y$.
\end{lemma}

\proof
We simply take $\beta$ to be the class in $H_1(\Sigma,\del\Sigma;\ZZ/2)$ of an arc meeting the closure of the descending manifold of the saddle point corresponding to an $x_i$ appearing with coefficient one once and all others zero times. See Figure \ref{fig4} for a picture of the descending manifolds.
\qed

Thus the ``cup-length'' of $HF_*(h,\ZZ/2;\eta)$ as an $H^1(S)$-module is two.

We will be taking a limit in which degenerate fixed points are allowed to appear. We have need of a Gromov compactness result appropriate for such a situation. We use Taubes's currents-in-the-target version of Gromov compactness (valid in dimension four) \cite[Proposition 3.3]{ta} as applied to compact subsets $[a,b] \times M_\phi$ in \cite[Lemma 9.9]{h}.

\begin{prop}[{\cite[Proposition 3.3]{ta},\cite[Lemma 9.9]{h}}]
\label{compact}
Let $u_k: \RR\times [0,1] \rightarrow \RR\times M_\phi$ be a sequence of holomorphic sections with energies bounded by some $E_0$. Then we can pass to a subsequence such that:
\be
\ii The $u_k$'s converge weakly as currents in $\RR\times M_\phi$ to a proper pseudoholomorphic map $u: \RR\times [0,1] \rightarrow \RR\times M_\phi$.
\ii For any compact $K \subset \RR \times M_\phi$, $$\quad\ \  \lim_{k\rightarrow\infty} \left[\mathrm{sup}_{x\in \Im(u_k)\cap K} \mathrm{dist}\left(x,\Im(u)\right) + \mathrm{sup}_{x\in \Im(u)\cap K} \mathrm{dist}\left(x,\Im(u_k)\right)\right] = 0. \quad\ \  \eop$$
\ee
\end{prop}

From this we see that we have $C^0$-convergence to the orbit corresponding to a fixed point at each end.

\begin{thm}
There must be at least two fixed points in the Nielsen class $\eta$ corresponding to such a fixed component $S$ even if we allow degenerate fixed points.
\end{thm}

\proof
Suppose there is just one, necessarily degenerate, fixed point in Nielsen class $\eta$. Perturb by a small Hamiltonian flow such that all fixed points are nondegenerate. Because $HF_*(h,\ZZ/2;\eta)$ has rank one in even degree, we have at least one index-one fixed point which survives in homology; choose one such and call it $y$. By Lemma \ref{beta}, for any index negative one fixed point $x$ which survives in homology, there is a $\beta \in H_1(S,\del S)$ such that $\del_\beta [x] = [y]$. Thus there are (at least) two flow-lines from $x$ to $y$, and there are two such that the difference between their classes, mapped to $H_1(S;\ZZ/2)$, is nonzero.

We now take a limit of small perturbations limiting to the degenerate situation. By Proposition \ref{compact}, after passing to a subsequence, the two flow-lines $u_k$ and $v_k$ above limit to holomorphic curves $u,v : \RR\times [0,1]\rightarrow \RR\times M_\phi$ which $C^0$-limit to the degenerate fixed point at each end. Thus each of these two limits give continuous loops in $\Gamma(M_\phi)_\eta$. We additionally see that $[u]-[v] = \lim_{k\rightarrow\infty} [u_k-v_k]$ in $H_1(\Gamma(M_\phi)_\eta \cong H_1(S)$. This latter limit is nonzero, and in particular at least one of $u, v$ is nonconstant and thus has positive energy. However, $\omega_\phi: H_1(\Gamma(M_\phi)_\eta) \rightarrow \RR$ is the zero map, so we see that each of $u, v$ has zero energy, a contradiction.
\qed

Summing over Nielsen classes, we conclude:

\begin{thm}
\label{degboundsection}
The minimum number of fixed points of an area-preserving map $\phi$ in a mapping class $h$ is given by:
$$ \left\{ \begin{array}{cc}
\#\left\{\eta:\ind(\eta)\neq 0\right\}  & h\ \textrm{is periodic or pseudo-Anosov} \\
\#\left\{\eta:\ind(\eta)\neq 0\right\} + A + B & h\ \textrm{is reducible} \end{array} \right.$$
where $A$ is as before and $B$ is the number of fixed annuli.
\end{thm}

\end{document}